\documentclass[12pt]{amsart}

\usepackage{amssymb}
\usepackage{bm}
\usepackage{graphicx}
\usepackage{enumerate}
\usepackage{url}
\newcommand{\excise}[1]{}%{$\star$\textsc{#1}$\star$}

\newtheorem{thm}{Theorem}[section]
\newtheorem{lemma}[thm]{Lemma}

\newtheorem{cor}[thm]{Corollary}

\theoremstyle{definition}

\newtheorem{example}[thm]{Example}
\newtheorem{remark}[thm]{Remark}
\newtheorem{defn}[thm]{Definition}

\numberwithin{equation}{section}

\DeclareMathOperator \mex{mex}
\DeclareMathOperator \lcm{lcm} % lcm
\newcommand \N {\mathbb{N}_0}
\newcommand \No {\mathbb{N}_1}
\newcommand \Z {\mathbb{Z}}
\newcommand \Pt {\tilde{P}}
\newcommand \po {\overline{p}}
\newcommand \Y {\mathcal{Y}}
\newcommand \Mo {\overline{M}}
\newcommand \Mu {\underline{M}}
\newcommand \C {\mathcal{C}}
\newcommand \D {\mathcal{D}}
\newcommand \Ct {\mathcal{\tilde{C}}}
\newcommand \Ch {\mathcal{\hat{C}}}
\newcommand \xt {\tilde{x}}
\newcommand \St {\tilde{S}}
\newcommand \Tt {\tilde{T}}
\newcommand \dt {\tilde{d}}
\newcommand \Sh {\hat{S}}
\newcommand \Th {\hat{T}}
\newcommand \dhat {\hat{d}}
\newcommand \B {\mathcal{B}}
\newcommand \gta {\Gamma} %{\widetilde{\Gamma}}
\newcommand \Et {E} %{\widetilde{E}}
\newcommand \Vt {V} %{\widetilde{V}}
\newcommand \Li {\mathrm{Li}}
\newcommand \Pl {\boldsymbol{+}}
\newcommand \Mi {\boldsymbol{-}}
\newcommand \Ro {\overline{R}}
\newcommand \pio {\overline{\pi}}
\newcommand \Kh {\hat{K}}
\newcommand \fw { \: | \: }
\newcommand \eentry {\underline{\phantom{+2}}}
\newcommand \scaletable[1] { \scalebox{0.8}{\ensuremath{ #1  } } }

\begin{document}%%%%%%%%%%%%%%%%%%%%%%%%%%%%%%%%%%%%%%%%%%%%%%%%%%%%%%%%
%%%%%%%%%%%%%%%%%%%%%%%%%%%%%%%%%%%%%%%%%%%%%%%%%%%%%%%%%%%%%%%%%%%%%%%%

\mbox{}

\title[On the additive period length of the S-G function]{On the additive period length of the Sprague-Grundy function of certain Nim-like games}

\author[J.~Askgaard]{Jens Askgaard}
\address{Denmark}
\email{jaskgaard\char`_math@outlook.dk}

\date{\today}

\begin{abstract}
We examine the structure of the additive period of the Sprague-Grundy function of Nim-like games, among them Wythoff's Game, and deduce a bound for the length of the period and preperiod.
\end{abstract}

\maketitle

% \setcounter{tocdepth}{1}
% \tableofcontents

\section{Introduction} %%%%%%%%%%%%%%%%%%%%%%%%%%%%%%%%%%%%%%%%%%%%%%%

Consider a sequence $\left( \Y_x \right)_{x=0}^\infty$ of finite subsets of $\Z$. They define a function $G: \N \to \N$ by 

$$ G(x) = \mex \left(  \{ G(x') | x' < x \} \cup Y_x \right) $$
where $\mex ( \Y ) = \min \left( \N \setminus \Y \right)$ is the minimal excluded operator, defined for any finite subset $\Y$ of $\N$ or $\Z$. Such functions appear as Sprague-Grundy functions in the study of certain impartial combinatorial games, such as Nim \cite{WinningWays}. We will impose the condition of additive periodicity on $\left( \Y_x \right)_{x=0}^\infty$, as defined below.

\begin{defn}[\textbf{Additive periodicity}] \label{d:AddPer}
Let $\left( \Y_x \right)_{x=0}^\infty$ be a sequence of finite subsets of $\Z$. Then $\left( \Y_x \right)_{x=0}^\infty$  is \emph{additively periodic}, if there exists $P' \in \N$ and $p \in \No$ such that for all $x \geq P'$, we have $\Y_{x + p} = \Y_x + p$. The uniquely determined smallest numbers $P', p$ for which this condition holds are called the \emph{preperiod length} and the \emph{period length} of $\left( \Y_x \right)_{x=0}^\infty$.\\
A function $G: \N \to \N$ is additively periodic, if there exists $\Pt \in \N$ and $\po \in \No$ such that for all $x \geq \Pt$, we have $G( x + \po) = G( x ) + \po$. Again, the smallest numbers $\Pt, \po$ for which this condition holds are called the preperiod length and the period length of $G$.
\end{defn}

Originally shown in \cite{Pink}, additive periodicity of $\left( \Y_x \right)_{x=0}^\infty$ implies that $G$ will be additively periodic as well. The motivation for this paper was to find the optimal bound for the period length of $G$, given only the period length of $\left( \Y_x \right)_{x=0}^\infty$, and the upper and lower bounds of the elements in the sequence $\left( Y_x - x \right)_{x=0}^\infty$.

To cover some practical cases from game theory where a finite number of the values of the Sprague-Grundy function are not defined by $ \left( \Y_x \right)$, we will introduce a seed in the definition of $G$, inspired by a similar definition in \cite{Abrams}.

\begin{defn}[\textbf{Nim sequence}] \label{d:NimSeq}
Let $\left( \Y_x \right)_{x=0}^\infty$ be a sequence of additively periodic finite subsets of $\Z$. Let $L \in \N$, and let $[g_0 , \ldots , g_{L-1}]$ be an $L$-tuple, where for all $x,x' \in \{ 0, \ldots , L-1 \}$, we have $g_x \in \N$, and $g_x = g_{x'}$ implies $x=x'$. Let $G: \N \to \N$ be defined by

$$G(x) = \left\{ \begin{array}{ll} 
g_x & \textrm{for $x<L$} \\
\mex( \{ G(x') | x' < x \} \cup \Y_x ) & \textrm{for $x \geq L$} .
\end{array} \right. $$
Then we call $G$ for a \emph{Nim sequence over} $\left( \Y_x \right)_{x=0}^\infty$, and $[g_0 , \ldots , g_{L-1}]$ is its \emph{seed}.
\end{defn}

Note that it is inconsequential what the sets $\left( \Y_x \right)_{x=0}^{L-1}$ are. So in any given example, we may redefine $L$ as $\max(L, P')$ and pretend that the sequence $\left( \Y_x \right)$ is additively periodic for all indices. Thus, we will silently assume that $P' = 0$ in this paper.

We study Nim sequences using their difference functions, the definition of which is given in \cite{Larsson}.

\begin{defn}[\textbf{Difference function}]
Given a function $G: \N \to \N$, we define its \emph{difference function} $d: \N \to \Z$ by $d(x) = G(x) - x$. If $G$ is additively periodic, we define its \emph{difference period} starting at $x \geq \Pt$ as the $\po$-tuple $[d(x), \ldots , d(x + \po - 1 )]$.

Given $\Y_x \subseteq \Z$, we use the notation $d( \Y_x) = \Y_x - x$.
\end{defn}

It is clear that $G$ is additively periodic, if and only if its difference function is periodic, so $d(x+\po) = d(x)$ for all $x \geq \Pt$. Similarly, if $\left( \Y_x \right)$ is additively periodic, then $d( \Y_{x+p} ) = d( \Y_x )$, so $\max_{ x \in \N} d( \Y_x ), \min_{ x \in \N} d( \Y_x )$ will be well-defined.

\begin{example} \label{e:Simple}
Let $\Mo \in \No$ and $\Mu \in -\No$. For all $x \in \N$, define 

$$\Y_x = \left\{ \Mu + 1, \ldots , \Mo - 1 \right\} + x ,$$
so $p=1$ and $\max d( \Y_x) = \Mo - 1$, $\min d(\Y_x) = \Mu + 1$. Let the seed be empty.

Then $\po = \Mo + |\Mu |$ with $\Pt = 0$, as the difference period becomes:

$$ [ \; \underbrace{\Mo, \Mo, \ldots, \Mo }_{|\Mu|  \; \textrm{times}  }, \underbrace{\Mu, \Mu, \ldots, \Mu }_{\Mo \; \textrm{times}  } \; ] .$$
\end{example}

Inspired by this example, we formally define $\Mo$ and $\Mu$ as constants determined by $\left( \Y_x \right)_{x=0}^\infty$.

\begin{defn}[\textbf{Difference bounds}]
Let $\left( \Y_x \right)_{x=0}^\infty$ be a sequence of additively periodic finite subsets of $\Z$. Then, we define the constants $\Mo = \max_{ x \in \N} d( \Y_x) + 1$, $\Mu = \min_{ x \in \N} d(\Y_x) - 1$, and $M = |\Mu| + \Mo$.
\end{defn}

The reason why we offset these maximal and minimal values with $1$ is that they will become the bounds for the difference function $d$ of $G$, as we will prove in this section. 
Now we can formulate the main result of this paper:

\begin{thm} \label{t:Main}
Let $G$ be a Nim sequence over $\left( \Y_x \right)_{x=0}^\infty$ that has period length $p$. Then $G$ is additively periodic, and the length $\po$ of its period is bounded by:
$$ \po \leq K_{\Mu , \Mo } \: p  $$
where $K_{\Mu , \Mo } \in \No$ is a constant, approximately defined by
$$ K_{\Mu , \Mo } \sim \exp \left( \sqrt{ \Li^{-1} \left( |\Mu | + \Mo - 1 \right) } \right) .$$
\end{thm}
Here, $\Li (x) = \int_2^x \frac{1}{\log t} \mathrm{d}t$ is the logarithmic integral, and the symbol $\sim$ signifies that $K_{\Mu , \Mo }$ has the same divergence speed as the expression.

The first part of this theorem, that $G$ is additively periodic, was proven by Pink in his diploma thesis. His proof was simplified by Dress and Flammenkamp in \cite{Pink}, which was followed by an even simpler proof by Landman in \cite{Landman}. These proofs all rely on pigeon-hole methods, which lead to larger bounds for $\po$. 

The canonical example of additive periodicity in game theory is Wythoff's Game, which we will briefly discuss.

\begin{example}[\textbf{Wythoff's Game}] \label{e:Wythoff}
For all $y \in \N$, define a Nim sequence $G_y$ by

\begin{eqnarray} G_y ( x ) & = & \mex \left( \left\{ G_y (x') | x' < x \right\} \cup \left\{ G_{y'} (x) | y' < y \right\} \right. \nonumber \\ 
& & {}\cup \left. \left\{ G_{ y-k } (x-k) | k \in \No , y-k \geq 0, x-k \geq 0 \right\}  \right) \label{eqn:Wythoff} .
\end{eqnarray}
Then $G_y (x)$ is the Sprague-Grundy value of the position $(x,y)$ in Wythoff's Game, as defined in \cite{WinningWays} and \cite{Landman}.

We give a matrix showing the first values of $G_y (x)$.

$$ \scaletable{ \begin{array}{c|ccccc ccccc cccc}
 (x,y) & 0 & 1 & 2 & 3 & 4 & 5 & 6 & 7 & 8 & 9 & 10 & 11 & 12 & 13 \\
\hline
0   & 0 & 1 & 2 & 3 & 4 & 5 & 6 & 7 & 8 & 9 & 10 & 11 & 12 & 13 \\
1   & 1 & 2 & 0 & 4 & 5 & 3 & 7 & 8 & 6 & 10 & 11 & 9 & 13 & 14\\
2   & 2 & 0 & 1 & 5 & 3 & 4 & 8 & 6 & 7 & 11& 9 & 10 & 14 & 12 \\
3   & 3 & 4 & 5 & 6 & 2 & 0 & 1 & 9 & 10 & 12 & 8 & 7 & 15 & 11 
\end{array} }$$
Wythoff showed that the positions, where $G_y (x) = 0$, are defined by 
$$(x,y) = \left( \lfloor k \varphi \rfloor , \lfloor k \varphi^2 \rfloor \right) \: \textrm{for all} \: k \in \N ,$$
where $\varphi$ is the golden ratio. That is, the zeroes approximately lie on the two diagonals emanating from the corner whose slopes equal $\varphi$ (see \cite{Wythoff}). However, the remaining values were considered to be chaotic (see \cite{WinningWays}), until the proof by Pink was published. Indeed, let $\Pt_y, \po_y$ be the preperiod length and the additive period length of $G_y$. Set 

$$ \Y_{x} = \left\{ G_{y'} (x ) | y' < y \right\} \cup \left\{ G_{ y-k } (x -k ) | k \in \No , y-k \geq 0, x-k \geq 0 \right\} .$$
With $ \displaystyle p = \prod_{y' < y} \po_{y'}$, assuming $\displaystyle x \geq \max_{y'<y} \Pt_{y'} + y$ we find $\Y_{x+p} = \Y_x + p$. We can now use expression (\ref{eqn:Wythoff}) plus Theorem \ref{t:Main} to create an induction proof which shows that each $G_y$ is additively periodic, as in \cite{Pink} or \cite{Landman}.

The following table shows the preperiod length and the period length for the first few rows.
$$ \scaletable{ \begin{array}{c|c|c|c}
y  & \Pt_y & \po_y & \textrm{Difference period} \\
\hline
 0  & 0 & 1 & [0] \\
 1  & 0 & 3 & [+1,+1,-2] \\
 2  & 0 & 3 & [+2,-1,-1] \\
 3  & 8 & 6 & [+2,+3,-2,-4,+3,-2] \\
\end{array} } $$
\end{example}

The lemmas in this introduction and their methods of proof are directly derived from Landman. As he uses a more narrow framework which only uses $d$ implicitly, we give the proofs in full; the reader familiar with any of the papers cited here are free to read only the lemmas in this section and skip their proofs.

\begin{lemma} \label{l:1}
Let $G$ be a Nim sequence over $\left( \Y_x \right)_{x=0}^\infty$. Then $G$ is a bijective function, meaning that $G$ is a permutation of $\N$.
\end{lemma}
\begin{proof}
Injectivity of $G$ follows from the definition of the $\mex$ operator and the requirement of the seed that $g_x \neq g_{x'}$ for $x \neq x'$.

Surjectivity of $G$: With $x, y \in \N$, there might be four causes why $G(x) \neq y$:
\begin{itemize}
\item[(i)] $G(x) < y$; this can at most occur $y$ times.
\item[(ii)] $y \in \Y_x$; this can at most occur $\sum_{j=0}^{p-1} \# \Y_j$ times.
\item[(iii)] With $x < L$, so $G(x)$ lies in the seed; this can at most occur $L$ times.
\item[(iv)] $G( x' ) = y$ for some $x'<x$. As (i), (ii) and (iii) occurs only finitely many times, (iv) must happen.
\end{itemize} 
\end{proof}

Thus, we can imagine $G$ as a greedy permutation that always chooses the smallest number not contained in $\Y_x$, as in \cite{Larsson}. 

\begin{lemma} \label{l:5}
Let $G$ be a Nim sequence over $\left( \Y_x \right)_{x=0}^\infty$. Then there exists some $C \in \N$ such that

$$ \max_{x \geq C} d(x) \leq \max( 0, \Mo) .$$
\end{lemma}

\begin{proof}
Set $\hat{M} = \max( 0, \Mo)$ First, assume $d(x) \leq \hat{M}$ for all indices $x < L$ in the seed. Assume there exists $x \in \N$ with $d(x) > \hat{M}$. As $\hat{M} \notin d( \Y_x)$ and $x + \hat{M} \in \N$, there must exist $x_0 < x$ such that $G( x_0 ) =x + \hat{M}$ due to the greedy definition of $G$. Set $a = \min \left\{ a' \in \N \fw G^{-1} (x+ \hat{M} - a' ) > x \right\} > 0$. For all $a' < a$, set $x_{a'} = G^{-1} (x+ \hat{M} - a' ) < x$.

We must have $x + \hat{M} - a \in \Y_{x_0} \cup \ldots \cup \Y_{x_{a-1} }$, and $x_{a-1} \leq x-a$. Then

\begin{eqnarray*}
\max_{x \in \N } ( d ( \Y_x) ) & < & \hat{M} \\
& \leq & x + \hat{M} - a - x_{a-1} \\
& \in & d( \Y_{x_{a-1} } )
\end{eqnarray*}
which is a contradiction.

If $d(x) > \hat{M}$ for some $x < L$, we first sort the seed so $g_0 < g_1 < \cdots < g_{L-1}$; as the $\mex$-operator depends on sets, this does not change the other values of $G$. 

If $G(L) > g_{L-1}$, we must have $G(L) = g_{L-1} + 1$ due to the bounds on $\left( Y_x \right)$. If $g(L) < g_{L-1}$, we have $d(L) < d(L-1) - 1$. In both cases, we add $G(L)$ to the seed, offset $L$ with $1$ and sort the seed again. Then, either $d(L-1)$ will remain the same, or $d(L-1)$ will decrease with $1$. We continue this process of adding more elements to the seed, until $d(L-1) = \hat{M}$, which must happen because $G$ is surjective. Now set $C = L$, and continue as in the first part of the proof.
\end{proof}

\begin{lemma} \label{l:6}
Let $G$ be a Nim sequence over $\left( \Y_x \right)_{x=0}^\infty$. Then there exists some $C \in \N$ with $C \leq \max_{x<L} g_x  + 1$ such that

$$ \min_{x \geq C} d(x) \geq \min( 0, \Mu) .$$
\end{lemma}
\begin{proof}
Set $\hat{M} = \min( 0, \Mu)$. Assume that $x \in \N$ is the first index where $d(x) < \hat{M}$. As in the proof of Lemma \ref{l:5}, we can assume that the seed is sorted by size, so $x \geq L$. We have assumed $G( x - 1) \geq x - 1 + \hat{M}$, so $G(x -1 ) > G(x)$. If $x > L$, then $G(x) \in \Y_{x-1}$, and

\begin{eqnarray*}
\hat{M} & > & d(x) + 1 \\
& = & x + d(x) - (x-1) \\
& \in & d( \Y_{x-1} )
\end{eqnarray*}
which is a contradiction.

If $x=L$, we add $G(L)$ to the seed, offset $L$ with $1$ and sort the seed again. Then, $d(L - 1)$ will decrease with $1$. We continue this process until $G(L - 1) < G(L)$, when we set $C=L$ and continue as in the first part of this proof.

This process could at most be repeated $\max_{x<L} g_x  - L$ times, where $L$ is the original size of the seed. So the maximal size of the increased seed is $L + \max_{x<L} g_x  - L = C-1$.
\end{proof}

\begin{cor} \label{c:6x}
Let $G$ be a Nim sequence over $\left( \Y_x \right)_{x=0}^\infty$. If $G$ is additively periodic with preperiod length $\Pt$, we have for all $x \geq \Pt$

$$ \min( 0, \Mu) \leq d(x) \leq \max( 0, \Mo) .$$
\end{cor}
\begin{proof}
Follows, as the difference values in the period repeat themselves.
\end{proof}

As we want $d$ to be bounded by a negative $\Mu$ and a positive $\Mo$, we need to deal with the degenerate case when it is not.

\begin{thm}
Let $G$ be a Nim sequence over $\left( \Y_x \right)_{x=0}^\infty$. If $\Mu \geq 0$ or $\Mo \leq 0$, then $G$ is additively periodic with period length $\po = 1$.
\end{thm}
\begin{proof}
Assume $\Mo \leq 0$. By Lemma \ref{l:5} we find $C \in \N$, so for all $x \geq C$ we have $G(x) \leq x$. As in the proofs of Lemma \ref{l:5}+\ref{l:6}, we keep extending the seed by adding $G(x)$ and sorting it, until we reach an index $L$ with $G(L) = L$. Then the seed will be equal to $\{ 0, \ldots, L-1 \}$, and $G(x) = x$ for all $x \geq L$, so the difference period becomes $[0]$.

The proof for $\Mu \geq 0$ is similar.
\end{proof}

\begin{remark} \label{rem:Main}
Knowing the bounds for the difference function, we can establish that $G$ is additively periodic. To generalize Landman's argument, by Lemma \ref{l:5}+\ref{l:6} and the surjectivity of $G$, we have for $x \in \N$ sufficiently large

$$\{ G( x' ) | x' < x - M \}  \subseteq \{ 0, \ldots , x - 1 + \Mu \} \subseteq \{ G(x') | x' < x \} .$$
Then, we calculate $G(x)$ as

\begin{eqnarray}
G(x) & = & \mex \big( \{ 0, \ldots , x - 1 + \Mu \} \nonumber \\
& & \cup \left( \{ G(x') | x -  M \leq x' < x \} \cap \{ x + \Mu, \ldots, x + \Mo - 1 \} \right) \cup \Y_x \big) \nonumber \\
& = & \mex \Big( \big( \{  -x, \ldots , - 1 + \Mu \} \label{e:Landman} \\ 
& & \cup \left( \{ G(x') - x | x -  M \leq x' < x \} \cap \{ \Mu, \ldots, \Mo - 1 \} \right) \cup d \left( \Y_x \right) \big)  +  x \Big) . \nonumber 
\end{eqnarray}
Now there are $2^M$ possibilities for what elements the set $\{ G(x') - x | x -  M \leq x' < x \} \cap \{ \Mu, \ldots, \Mo - 1 \}$ can contain. Using the pigeon-hole principle, we find $r , r' \in \N$, $r < r'$ such that 

\begin{eqnarray*}
& & \{ G(x') -  (x + rp) | x + rp -  M \leq x' < x + rp \} \cap \{ \Mu, \ldots, \Mo - 1 \} \\
& = & \{ G(x') - (x + r'p) | x + r'p -  M \leq x' < x + r'p \} \cap \{ \Mu, \ldots, \Mo - 1 \} .
\end{eqnarray*}
As $d \left( \Y_{x + rp} \right) = d \left( \Y_{x + r'p} \right)$, we use (\ref{e:Landman}) to show that $G( x + rp ) = G( x + r'p )$. Then, we continue inductively to show that for all $y \geq x + rp$, we have $G( y ) = G( y + (r' - r)p )$.
\end{remark}

Landman's proof implies that $K_{\Mu , \Mo } \leq 2^M $. In the next section, we will deduce another proof of the additive periodicity which, while still using pigeon-hole methods, yields a smaller bound.

\section{Periodicity conditions} %%%%%%%%%%%%%%%%%%%%%%%%%%%%%%%%%%%%%%%%%%%%%%%%%%%%%%%%%%%%

Throughout this section, $\left( \Y_x \right)_{x=0}^\infty$ stands for a sequence of additively periodic finite subsets of $\Z$ with period length $p$ that defines the difference bounds $\Mu \in -\No, \Mo \in \No$. Also, $G$ is a Nim sequence over $(\Y_x)$ with a seed of length $L$. We mention that while $p$ always stands for the smallest number that satisfies $\Y_x + p = \Y_{x+p}$, if we find an integer $p' \in \No$ such that $G(x + p') = G(x) + p'$ for all $x \geq \Pt'$ for some $\Pt' \in \N$, then $p', \Pt'$ may not be the true period / preperiod length. However, the true period length $\po$ of $G$ will divide $p'$.

We begin with a simple lemma that contains an important definition.
\begin{lemma}[\textbf{Exclusion Lemma}]
Let $(x,y)$ be an inversion of $G$; that is, $x<y$ with $G(x)>G(y)$. Then for all $a \in \Z$ with $ x + ap \geq L$, we have
$$ d( x + ap) \neq G(y) - x.$$
We say that (the difference value) $G(y)-x$ is \emph{excluded} at the index $x$ (mod $p$).
\end{lemma}
\begin{proof}
As $G$ is a greedy permutation, we must have $G(y) \in \Y_x$. Then $d(x + ap) \notin d( Y_{x+ ap} )= d( Y_x ) \ni G(y) - x$.
\end{proof}

We mention two simple facts. First, the excluded value is bounded by $\Mu < G(y) - x < \Mo$. Second, the excluded value can be calculated without knowing the exact values of $y$ and $G(x)$.

\begin{lemma} \label{l:4a}
Let $R \in \No$. Then $G$ is additively periodic with period length $\po$ that divides $Rp$, if and only if there exists $x \in \N$ that fulfills these two conditions:
\begin{itemize}
\item[a)] $\left\{G( x + k ) \bmod Rp \fw 0 < k < Rp \right\} = \left\{  0, \ldots , Rp-1 \right\} $,
\item[b)] $\left\{ G(x') | x' < x  \right\} = \bigcup_{k=0}^{Rp - 1} U_k \cap \N$
\end{itemize}
where $U_k = \left\{ G(x + k ) - aRp \fw a \in \No \right\} $.
\end{lemma}
\begin{proof}
Assume a) and b) are true. For $k, k' \in \{ 0, \ldots, Rp - 1 \}$, as $G( x + k ) \not\equiv G( x + k') \pmod {Rp}$ when $k \neq k'$, the sets $\left( U_k \right)_{k=0}^{Rp-1}$ are mutually disjoint. Then:

\begin{eqnarray*}
G( x + Rp ) & = & \mex \left( \bigcup_{k=0}^{Rp-1} \big( U_k \cup \left\{ G( x + k ) \right\} \big) \cup \Y_{x+Rp} \right) \\
& = & \mex \left( \bigcup_{k=0}^{Rp-1} \big( \left\{ G(x + k ) - aRp \fw a \in \N \right\} \big)  \cup \left( \Y_{x} + Rp \right) \right) \\
& = & \mex \left( \left( \bigcup_{k=0}^{Rp-1}U_k \cup \Y_x  \right) + Rp \right) \\
& = & G( x ) + Rp .
\end{eqnarray*}

For $x+1$, a) holds because:

\begin{eqnarray*}
& &\left\{G( x + 1 + k ) \bmod Rp \fw 0 < k < Rp \right\} \\ 
& = & \left( \left\{G( x + k ) \bmod Rp \fw 0 < k < Rp \right\} \setminus \left\{ G( x )  \bmod Rp \right\} \right) \\ 
& & {}\cup \left\{ G( x + Rp ) \bmod Rp \right\} \\
& = & \left\{G( x + k ) \bmod Rp \fw 0 < k < Rp \right\} .
\end{eqnarray*}
Then $\left\{ G(x + 1 + k ) - aRp \fw a \in \No \right\}  = U_{k+1}$ for $k < Rp-1$, and \\ $\left\{ G(x + 1 + Rp - 1) - aRp \fw a \in \No \right\}  = \{ G( x )\} \cup U_0$, so b) holds as well. We can now continue inductively to show that for all $x' \geq x$, we have $G( x' + Rp ) = G( x') + Rp$, so $x$ must be (larger than) the preperiod length.

Assume now that $G$ is additively periodic. Set $x = \Pt$. If a) does not hold, there must exist $k, k' \in \{ 0, \ldots, Rp - 1\}$ with $k \neq k'$ and $G(x + k) \equiv G( x + k') \pmod {Rp}$; assuming $G( x + k ) < G( x + k')$, we have $G( x + k ) = G( x + k' ) + CRp = G( x + k' + CRp)$ for some $C \in \No$. As $|k - k'| < Rp$, this is a contradiction.

If b) does not hold, there must exist $x' < x, C \in \N$ such that $G(x' ) = G( x +k + CRp ) = G( x + k ) + CRp$, which is also a contradiction.
\end{proof}

We now define for any $ x \in \N$ the two sets of indices where the permutation $G$ ascends above $x$, and where it descends below $x$.

\begin{defn}[\textbf{Cut}] \label{d:Cut}
Let $x \in \N$. We define the \emph{cut after $x$} as the pair of sets $(S_x , T_x)$, where

\begin{eqnarray*}
 S_x & = & \left\{ x' \in  \N \fw x' < x + \frac{1}{2} , G(x') > x + \frac{1}{2} \right\}  , \\
 T_x & = & \left\{ x' \in  \N \fw x' > x + \frac{1}{2} , G(x') < x + \frac{1}{2} \right\} .
\end{eqnarray*}
With each cut, we associate the sets $S_x^* = G( T_x)$ and $T_x^* = G( S_x )$.

For $x, y \in \N$, we call the sets $S_x, S_y$ \emph{equivalent}, written $S_x \sim S_y$, if and only if $d( S_x ) = d( S_y)$.
\end{defn}

It is clear that $\sim$ is an equivalence relation, and we will use it to compare any two sets with an index in $\N$ or $\Z$.

\begin{remark}
Since $G$ is bijective, we always have

$$\# S_x = \# T_x = \# S_x^* = \# T_x^* < \infty .$$
\end{remark}

\begin{lemma} \label{l:13}
Let $R \in \No$. Then $G$ is additively periodic with period length $\po$ that divides $Rp$, if and only if there exists $x \geq L$ such that $S_x^* \sim S_{x+Rp}^*$, and $T_x^* \sim T_{x+Rp}^*$.
\end{lemma}
\begin{proof}
Assume the two sets are equivalent for some $x \geq L$, then we will show the conditions from Lemma \ref{l:4a} hold for $x+1$. We also set $R=1$; for $R>1$, simply replace $p$ with $Rp$ in the following argument. Assume there exist $x_0' , x_0'' \in \left\{ x+1, \ldots, x + p \right\}$, $x_0' \neq x_0''$ with $G( x_0' ) \equiv G( x_0'' ) \pmod p $, $G( x_0' ) < G( x_0'' )$.

If $G( x_0' ) > x + p$, that means $x_0' \in S_{x+p}$, by assumption there exists $x_1' \in S_x$ with $G( x_1' ) = G( x_0' ) - p$. If $ x_1' \in S_{x+p}$, we find $x_2' \in S_x$ with $G( x_2' ) = G( x_0' ) - 2p$. We continue inductively, until we find $a \in \N$, $x_a' \in S_x$ with $G( x_0' ) - ap = G( x_a') \in \left\{ x+1, \ldots, x+p \right\}$.

Correspondingly, if $x_0' \in T_{x}$, we find $x_1' \in T_{x+p}$ with $G(x_1') = G(x') + p$, and we continue inductively, until we find $a \in \N$, $x_a' \in T_{x+p} $ with $G( x_0' ) + ap = G( x_a') \in \left\{ x+1, \ldots, x+p \right\}$.

We use the same process with $x_0''$ to find $b \in \Z$ such that $G( x_0'' ) + bp = G( x_b'') \in \left\{ x+1, \ldots, x+p \right\}$, so $G( x_a' ) \equiv G( x_b'') \pmod p$. If $x_0 ' \in T_x $, $x_0'' \in S_{x+p}$, then $G( x_a') = G(x_b'')$ as both $G$-values are bounded between $x+1$ and $x+p$, . As $x_a' \neq x_b''$, this contradicts the injectivity of $G$ . If instead $x_0', x_0'' \in S_{x+p}$, we find $x_b'' \in S_x$ with $G( x_0' ) = G( x_0'' ) - bp =  G( x_b'') $, again contradicting the injectivity of $G$. If $x_0', x_0'' \in T_x$, we find $x_a' \in T_{x+p}$ with $G( x_0'' ) = G( x_0' ) + ap =  G( x_a' ) $, again with a contradiction. Thus, condition a) is proven.

Now, let $x_0' \in \left\{ x+1, \ldots, x+p \right\}$, $a \in \Z$, $x_0'' \in \N$ such that $G( x_0'' ) = G( x_0' ) - ap$.  If $ a < 0$, we must have $x' < x''$ for condition b) to hold. If we assume $x' > x''$, we have $x_0' \in T_x$ or $x_0'' \in S_{x+p}$ (possibly both), and as above, we can find a contradiction by showing $G$ is not injective. If $a > 0$, we must have $x_0' > x_0''$ for b) to hold. If we assume $x_0' < x_0''$, we have $x_0' \in S_{x+p}$ or $x_0'' \in T_{x+p}$, and the contradiction follows in the same manner.

Now assume $G$ is additively periodic. Then, for $x \geq \Pt + \Mo$, it is trivial to show the other implication.
\end{proof}

\begin{lemma} \label{l:13a}
Let $R \in \No$ and $x \geq L$. Then:

\begin{itemize}
\item[a)] If $S_x \sim S_{x + Rp}$, then $T_x^*  \sim T_{x + Rp}^*$.
\item[b)] If $S_x^* \sim S_{x + Rp}^*$, then $T_x \sim T_{x + Rp}$.
\end{itemize}
\end{lemma}
\begin{proof}
Set $R=1$; for $R>1$, replace $p$ with $Rp$ in the following argument.

To prove a), let $N = \# S_x$, and let $\left\{ x_1, \ldots, x_N \right\} = S_x$ be in ascending order. First assume $G( x_1 ) + p < G( x_1 + p )$, so $G( x_1 ) + p > x + p$. If $G(x_1) + p \in T_{x+p}^*$, its inverse must belong to $S_x \setminus \left\{ x_1 \right\}$. So $G^{-1}(G(x_1) + p) \in \{ x_2 + p, \ldots, x_N + p\} \cup \{ x + p + 1, \ldots \}$, which means we have an inversion $( x_1 + p, G^{-1}(G(x_1) + p) )$ that excludes $G(x_1 ) + p - (x_1 + p ) = d( x_1 )$ at the index $x_1 + p$. Subtract $p$ from the index to see that $d(x_1)$ is excluded at its own index, which is a contradiction. If $G( x_1 ) + p > G( x_1 + p )$, we find an inversion $(x_1, G^{-1}(G( x_1 + p ) - p ) )$ in the same manner that would exclude $d(x_1 + p )$ at the index $x_1$, also a contradiction.

As we now have $G( x_1 ) + p = G( x_1 + p )$, let us look at $x_2$. If $G( x_2 ) + p < G( x_2 + p )$, we would have $G^{-1}(G(x_2) + p) \in \{ x_3 + p, \ldots, x_N + p\} \cup \{ x + p + 1, \ldots \}$, and we would find an impossible exclusion. The same thing happens if $G( x_2 ) + p > G( x_2 + p )$. So $G( x_2 ) + p = G( x_2 + p )$, and we continue inductively for each $n \in \left\{ 3, \ldots, N \right\}$ to show $G( x_n ) + p = G( x_n + p )$, which implies $T_x^*  \sim T_{x + Rp}^*$.

To prove b), let $\left\{ y_1, \ldots, y_N \right\} = S_x^*$ be in ascending order and continue using similar ideas. For instance, if $G^{-1}( y_1 ) + p < G^{-1}( y_1 + p )$, then $( G^{-1}( y_1 ) + p , G^{-1}( y_1 + p ) )$ is an inversion that would exclude $d( G^{-1} ( y_1 ) )$ at the index $G^{-1} ( y_1 ) + p$.
\end{proof}

Note that the reverse implications of a) and b) are not true in general.

The proof technique of the previous lemma will be re-used in the first part of the following proof.

\begin{lemma} \label{l:13c}
Let $R \in \No$. Then $G$ is additively periodic with period length $\po$ that divides $Rp$, if and only if there exists $x \geq L$ such that $T_x \sim T_{x+Rp}$, and $T_x^* \sim T_{x+Rp}^*$.
\end{lemma}
\begin{proof}
(You may replace $p$ with $Rp$ in the following argument). Assume the sets are equivalent. If $T_x \neq \emptyset$, let $\hat{x} = \max T_x$, else let $\hat{x} = x$. We show that for all $x_n \in \left\{ x_1, \ldots, x_N \right\} = \left\{ x+1, \ldots, \hat{x} - 1 \right\} \setminus T_x$, we have $G( x_n ) + p = G( x_n )$, where $x_1 , \ldots , x_N$ are in ascending order, and $N \in \N$ is the length of the sequence.

First assume $G( x_1 ) + p < G( x_1 + p)$. As $x_1 \notin T_x$, we have $G( x_1 ) > x$, so $G( x_1 ) \notin T_x^*$. Per assumption we have $G( x_1  )+ p  \notin T_{x+p}^*$, so $G^{-1}(G(x_1) + p) \in \{ x_2 + p, \ldots, x_N + p\} \cup \{ \hat{x} + p , \hat{x} + p + 1, \ldots \}$. Then $( x_1 + p, G^{-1}(G(x_1) + p) )$ is an inversion that excludes $d( x_1 )$ at the index $x_1 + p$, which is a contradiction. If we assume $G( x_1 ) + p > G( x_1 + p)$, the inversion $(x_1, G^{-1}( G( x_1 + p ) - p ) )$ impossibly excludes $d(x_1 + p )$ at the index $x_1$.

We continue inductively for each $n \in \{ 2, \ldots, N \}$, and find $G( x_n ) + p = G(x_n + p )$. We now realize that $T_{\hat{x} }^* \sim T_{\hat{x} + p}^*$, and $S_{\hat{x} }^* \sim S_{\hat{x} + p}^*$ is shown by

\begin{eqnarray*}
& & S_{\hat{x} }^* + p \\
& = & \left( \left\{ x+1, \ldots, \hat{x} - 1 \right\} \setminus \left( T_x^* \cup \left\{ G( x_1 ), \ldots, G( x_N ) \right\} \right) \right) + p \\
& = & \left\{ x+1 + p, \ldots, \hat{x} - 1 + p \right\} \setminus \left( T_{x+p}^* \cup \left\{ G( x_1 + p), \ldots, G( x_N + p) \right\} \right) \\
& = & S_{ \hat{x} + p }^* .
\end{eqnarray*}
So the conditions of Lemma \ref{l:13} are fulfilled for $\hat{x}$.

If we assume $G$ is additively periodic, it is trivial to show the other implication for $x$ sufficiently large.
\end{proof}

\begin{cor} \label{c:13c'}
Let $R \in \No$. Then $G$ is additively periodic with period length $\po$ that divides $Rp$, if and only if there exists $x \geq L$ such that $S_x \sim S_{x + Rp}$, and one of the following statements holds:

\begin{itemize}
\item[a)] $T_x  \sim T_{x + Rp}$; or
\item[b)] $S_x^* \sim S_{x + Rp}^*$.
\end{itemize}
\end{cor}
\begin{proof}
Combine Lemma \ref{l:13a} with Lemma \ref{l:13c} for a), and with Lemma \ref{l:13} for b).
\end{proof}

We can now prove the additive periodicity of $G$ with an improved bound for the period length. The same result was given by Dress and Flammenkamp in \cite{Pink}, which used a somewhat different proof technique.

\begin{thm} \label{t:Main2}
Any Nim sequence $G$ is additively periodic, and the length $\po$ of its period is bounded by:

$$\po \leq \binom{M}{ \min \left( | \Mu |, \Mo \right) } p .$$
\end{thm}
\begin{proof}
Let $x \in \N$ be sufficiently large so the bounds from Lemma \ref{l:5} and \ref{l:6} hold. Then look at the sets $\left( d( S_{x + rp} ), d( T_{x+ rp} ) \right)_{r=0}^\infty $. By the pigeon-hole principle, there exists $r, r' \in \N$, $r < r'$ such that $d( S_{x + rp} ) = d( S_{x + r'p} )$ and $d( T_{x + rp} ) = d( T_{x + r'p} )$, when $G$ becomes additively periodic with period length less or equal to $(r' - r)p$ by Corollary \ref{c:13c'}.

The maximal value of $r' - r$ can be calculated by counting all possibilites for $( d( S_x ), d( T_x ) )$. We have $\# S_x \leq \Mo$ and $\# T_x \leq | \Mu |$. By taking a sum over $i = \# S_x = \# T_x$, we get

$$ r' - r \leq \sum_{i=0}^{ \min \left( | \Mu |, \Mo \right) } \binom{\Mo}{i} \binom{ |\Mu | }{i} = \binom{M}{ \min \left( | \Mu |, \Mo \right) } $$
using the Chu-Vandermonde identity.
\end{proof}

Note that we could replace all references to $( S_x, T_x)$ in this proof with any one of the pairs $( S_x^*, T_x), ( S_x, T_x^*), ( S_x^*, T_x^*)$ and still reach the same result.

\section{Cut sets} %%%%%%%%%%%%%%%%%%%%%%%%%%%%%%%%%%%%%%%%%%%%%%%%%%%%%%%%%%%%%%%%%%%%%%

To move beyond the pigeon-hole proofs, we need to create a framework where cuts exist without reference to the permutation $G$ or the sets $( \Y_x )$, so $x \in \Z$ becomes a dummy variable.

\begin{defn}[\textbf{Cut set}] \label{d:CutSet}
Let $x \in \Z$ and $S_x, T_x \subseteq \Z$, and let $d: S_x \cup T_x \to \Z$. Then we define $\C_x = ( S_x , T_x , d)$ as a \emph{cut set}, if and only if these four requirements hold:
\begin{itemize}
\item[a)] For all $y \in S_x$, we have $y < x + \frac{1}{2}$ and $y + d(y) > x + \frac{1}{2}$.
\item[b)] For all $y \in T_x$, we have $y > x + \frac{1}{2}$ and $y + d(y) < x + \frac{1}{2}$.
\item[c)] $\# S_x = \# T_x < \infty $.
\item[d)] For all $y, y' \in S_x \cup T_x$, we have $y + d(y) = y' + d(y')$, if and only if $y = y'$.
\end{itemize}

With each cut set, we associate the sets 
$$S_x^* = \left\{ y + d(y) | y \in T_x \right\}, \quad T_x^* = \left\{ y + d(y) | y \in S_x \right\} .$$

We say that two cut sets $\C_x$ and $\Ct_{\xt} = \left( \St_{\xt} , \Tt_{\xt} , \dt \right)$ are \emph{equivalent}, written $\C_x \sim \Ct_{\xt}$, if and only if $T_x \sim \Tt_{\xt}$, $T_x^* \sim \Tt_{\xt}^*$, and $d(y) = \dt( y - x + \xt )$ for all $y \in S_x \cup T_x$.
\end{defn}

It is clear that $\sim$ is an equivalence relation, and each element in an equivalence class $\left[ \C_0 \right]$ corresponds to an $x \in \Z$.

We now need to replace the sets $( \Y_x )$. For this, we will use exclusions.

\begin{defn}[\textbf{Exclusions and possible zeros}] \label{d:Exc}
Let $\C_x = ( S_x , T_x , d)$ be a cut set. Let $x' \in S_x$. For all $x'' \in S_x$ with $x' < x''$ and $x' + d( x' ) > x'' + d( x'' )$, the difference value $x'' + d( x'' ) - x'$ is \emph{excluded} at the index $x'$, and for all $y \in \left\{ x+1, \ldots , x' + d( x' ) - 1 \right\} \setminus T_x^*$, the value $y - x'$ is excluded at the index $x'$. These are the \emph{positive exclusions}.

Let $y'' \in T_x$. For all $y' \in \left\{ x+1, \ldots , y'' - 1 \right\}$ with either $y' \notin T_x$, or $y' \in T_x$ with $y' + d( y' ) > y'' + d( y'')$, the difference value $y'' + d( y'' ) - y'$ is excluded at the index $y'$. These are the \emph{negative exclusions}.

We say there is a \emph{possible zero} (in difference value) after the cut, if and only if $x+1 \notin T_x \cup T_x^*$. We say there is a possible zero before the cut, if and only if $x \notin S_x \cup S_x^*$.
\end{defn}

We cannot find exlusions of the value $0$ by looking at one cut set only. For this, we need to have two cut sets in succession.

\begin{defn}[\textbf{Direct successor}] \label{d:DirSuc}
Let $\C_x , \Ct_{\xt}$ be two cut sets. We call $\Ct_{\xt}$ for a \emph{direct successor} of $\C_x$, if and only if there exists a cut set $\C_{x+1} = \left( S_{x+1}, T_{x+1}, \dt \right)$, which fulfills $\C_{x+1} \sim \Ct_{\xt}$ and the following conditions:

Initially, we set $S_{x+1} = S_x \setminus \left\{ x' \in S_x \fw x' + d( x' ) = x+1 \right\}$, and $T_{x+1} = T_x \setminus \{ x+1 \}$. Then we adjoin one or zero indices to $S_{x+1}$ and/or $T_{x+1}$, based on these criteria:

\begin{itemize}
\item[(i)] With $x+1 \in T_x \cap T_x^*$, we cannot adjoin any indices.
\item[(ii)] With $x+1 \in T_x^* \setminus T_x$, we adjoin $x+1$ to $S_{x+1}$, and $\dt ( x+1 ) > 0$ can be set to any value where $x+1 + \dt (x+1) \notin T_x^*$. 
\item[(iii)] With $x+1 \in T_x \setminus T_x^*$, we adjoin $x' \in \{ y \in \Z \fw y > x+1 \} \setminus T_x $ to $T_{x+1}$, and set $\dt (x') = x' - (x+1)$.
\item[(iv)] With $x+1 \notin T_x \cup T_x^*$, meaning that $\C_x$ has a possible zero after the cut, we have two options:
\item[(ivA)] We adjoin $x+1$ to $S_{x+1}$ and a new index $x'$ to $T_{x+1}$ such that $\dt (x+1)$ fulfills the conditions from (iii), and $x'$ fulfills the condtions from (iv). Then, we say that the value $0$ is \emph{excluded} at the index $x+1$.
\item[(ivB)] We add nothing, meaning that $\C_{x+1}$ has a possible zero before the cut. Then, we say that $\left( \C_x, \C_{x+1} \right)$ has a \emph{matching zero}.
\end{itemize}
\end{defn}

Note that these conditions will guarantee that $\# S_{x+1} = \# T_{x+1}$.

To handle exclusions, we need to define a collection of several cut sets.

\begin{defn}[\textbf{Cut set with several rows}] \label{d:CutSet2}
Let $x \in \Z$ and $R \in \No$. Then we define $\C_x = \left( \C_{x,r} \right)_{r=0}^{R-1}$ as a \emph{cut set with $R$ rows}, if for all $r, r' \in \{ 0, \ldots , R-1 \}$, where $\C_{x,r} = \left( S_{x,r}, T_{x,r}, d_r \right)$ is a cut set, we have:

\begin{itemize}
\item[a)] $T_{x,r} \sim T_{x,r'}$ and $T_{x,r}^* \sim T_{x,r'}^*$, if and only if $r = r'$.
\item[b)] If the difference value $k \in \Z$ is excluded at the index $x'$ in $\C_{x,r}$, then $x' \in S_{x,r'} \cup T_{x,r'}$ implies that $d_{ r' } (x' ) \neq k$.
\end{itemize}

We say that two cut sets $\C_x$ and $\Ct_{\xt} = \left( \Ct_{\xt, r} \right)_{r=0}^{R-1}$ with the same number of rows are \emph{equivalent}, written $\C_x \sim \Ct_{\xt}$, if and only if $\C_{x,r} \sim \Ct_{ \xt , r}$ for all $r \in \{0, \ldots R-1 \}$.
\end{defn}

\begin{defn}[\textbf{Direct successor}] \label{d:DirSuc2}
Let $\C_x , \Ct_{\xt}$ be two cut sets with $R$ rows. We call $\Ct_{\xt}$ for a \emph{direct successor} of $\C_x$, if and only if there exists a cut set  $\C_{x+1} = \left( \C_{x+1, r} \right)_{r=0}^{R-1}$ with $C_{x+1, r} = \left( S_{x+1, r}, T_{x+1, r}, \dt_r \right)$, which fulfills $\C_{x+1} \sim \Ct_{\xt}$ and the following conditions for all $r, r' \in \{ 0, \ldots, R-1 \}$:

\begin{itemize}
\item[a)] $C_{x+1, r}$ is a direct successor of $C_{x,r}$.
\item[b)] If $\left( C_{x,r}, C_{x+1, r} \right)$ excludes the difference value $0$ at the index $x+1$, \\then $\left( C_{x,r'}, C_{x+1, r'} \right)$ cannot have a matching zero.
\end{itemize}
\end{defn}

\begin{defn}[\textbf{Path}]
Let $R \in \No$, $\Mu \in -\No$ and $\Mo \in \No$. Then $\gta ( R, \Mu, \Mo ) = \left( \Vt ( R, \Mu, \Mo ) , \Et ( R, \Mu, \Mo ) \right)$ is the digraph, where the vertices $\Vt ( R, \Mu, \Mo )$ are all cut sets $\C_x =\left( S_{x+1, r}, T_{x+1, r}, \dt_r \right)_{r=0}^{R-1}$ with $R$ rows which fulfill $\bigcup_{r=0}^{R-1} d_r \left( S_{x,r} \cup T_{x,r} \right) \subseteq \left\{ \Mu , \ldots , \Mo \right\}$, and the edges are defined by the direct successors from Definition \ref{d:DirSuc2}.

With $J \in \N$, we call $\C_0 , \C_1 , \ldots, \C_J $ for a \emph{path} in $\gta ( R, \Mu, \Mo )$, if and only if $\left( \C_j ,  \C_{j+1} \right) \in \Et ( R, \Mu, \Mo )$ for all $j \in \{ 0 , \ldots , J-1 \}$. We also use the terms \emph{ancestor} and \emph{successor}, defined in the standard way for a digraph.
\end{defn}

We will usually abbreviate the notation as $ \gta = \left( \Vt , \Et \right)$.

We think of a cut set as an $R \times M$ matrix of difference values. The entries with undefined difference functions are left blank. A path with $J$ elements can be considered as an $R \times (M+J-1)$ matrix that defines $R$ common difference functions $d_0 , \ldots , d_{R-1}$, one for each row.

\begin{example} \label{e:CutSet}
Here is an example with a cut set and a possible successor in $\gta( 3, -2, 3 )$.

$$\scaletable{ \begin{array}{ccc|cc}
+3 & \eentry & \eentry & -1 & \eentry \\
\eentry & \eentry & +2 & -2 & \eentry \\
\eentry & \eentry & \eentry & \eentry & \eentry
\end{array} 
\quad \hookrightarrow \quad 
\begin{array}{ccc|cc}
\eentry & \eentry & \eentry & \eentry & \eentry  \\
\eentry & +2 & \eentry  & -1 & \eentry\\
\eentry & \eentry & +2 &  -1 & \eentry 
\end{array} }$$
(The horizontal lines mark the cut between positive and negative difference values.)
\end{example}

\begin{example} \label{e:R22}
Let $R \in \No$, and let us look at the digraph $\gta( R, -2 ,2 )$. It follows from Theorem \ref{t:Main2} that if $R > \binom{4}{2} = 6$, then the digraph is empty. If $R=6$, we can construct an infinite path in $\gta$ with a cycle of four cut sets $\C_0 , \ldots, \C_3, \C_4 \sim \C_0$. Here is an example:

$$ \scaletable{ \begin{array}{cc|cccc|cc}
+2 & +2  & -2 & -2 & +2 & +2 & -2 & -2 \\
\eentry & +2 & +2 & -2 & -2 & +2 & \eentry & -2 \\
\eentry & \eentry & +2 & +2 & -2 & -2 & \eentry & \eentry \\
+2 & \eentry & -2 & +2 & +2 & -2 & -2 & \eentry \\
\eentry & +2 & -2 & +2 & -2 & +2 & -2 & \eentry \\
+2 & \eentry & +2 & -2 & +2 & -2 & \eentry & -2 \\
\end{array} }$$
(The first horizontal line marks the end of $S_{0,r}$, and the second line marks the beginning of $T_{4,r}$.)

It can be shown that these four cut sets form a \emph{connected component}, and that all connected components in $\gta( 6, -2, 2 )$ consist of exactly four cut sets, which may be created by applying the same permutation of the six rows to our four cut sets. Note that the path has two separate row cycles, one $4$-cycle which cycles the top four rows, and one $2$-cycle which cycles the lower two rows. This implies that a Nim sequence with $\Mu = -2$, $\Mo = 2$ and a difference period of length $6p$ cannot exist, as we will show.
\end{example}

\begin{defn}[\textbf{Cycled cut set}] \label{d:CycCut}
Let $\C_x , \Ct_{\xt}$ be two cut sets with $R$ rows. We call $\Ct_{\xt}$ for the \emph{cycled cut set of $\C_x$}, if and only if for all $r \in \{ 0, \ldots , R-1 \}$

$$ \C_{x,r} \sim \Ct_{\xt, \pi( r )} ,$$
where $\pi$ is the row permutation that cycles all rows of $\C_x$ upwards, defined by $\pi (r) = ( r-1 ) \bmod R$.
\end{defn}

We now establish the connection between this abstract setup and the Nim sequence $G$. We need the additive period length of $G$ to be a multiple of $p$, which explains the clumsy definition of $R$ in the following lemma.

\begin{lemma} \label{l:17}
Let $G$ be a Nim sequence over $(\Y_x)$ that has additive period length $p$. Let $G$ have difference function $d$, preperiod length $\Pt$ and additive period length $\po$, and let $(S_x , T_x)$ be the cut after $x \in \N$ defined by $G$. 

Set $R= \frac{\lcm(p, \po ) }{p}$. Then there exists a path of cut sets $\C_0, \ldots , \C_p$ with $R$ rows, where $\C_j = \left( S_{j,r} , T_{j,r} , d_r \right)_{r=0}^{R-1}$, which fulfills:

\begin{itemize}
\item[a)] $\C_p$ is the cycled cut set of $\C_0$;
\item[b)] For all $x \geq \Pt$, where $j \in \{ 0, \ldots, p-1 \}$, $r \in \{ 0, \ldots, R-1 \}$, $C \in \N$ is uniquely defined by $x = j + rp + CRp$, we have:

$$ S_{j,r} \sim S_x \: , \: T_{j,r} \sim T_x \: ;$$
\item[c)] With $k \in \Z$, where $j+k \in S_{j,r} \cup T_{j,r}$, we have $d_r (j+k) = d( x+k )$.
\end{itemize}
\end{lemma}
\begin{proof}
Set $x' = \left( \left\lceil \frac{\Pt}{p} \right\rceil + C' \right) p$ for some sufficiently large $C' \in \N$. Write the difference period of $G$, as it begins at the index $x'$, in a $R \times p$ matrix. Extend the matrix to the left with the elements of $d( S_{x'+rp} )$ and to the right with the elements of $d( T_{x'+p-1+rp} )$ in each row with index $r$, as in Example \ref{e:R22}, so the path $\C_0 , \ldots, \C_{p-1}$ will be represented by the matrix. Conditions a) and b) now follow directly.
\end{proof}

\begin{lemma} \label{l:18}
Let $\C_0 , \ldots , \C_p$ be a path of cut sets with $R$ rows, where $\C_p$ is the cycled cut set of $\C_0$, and none of $\C_1, \ldots , \C_{p-1}$ is the cycled cut set of $\C_0$.

Then there exist an additive periodic sequence $(\Y_x)_{x=0}^\infty$ of finite subsets of $\Z$ with period length $p$, and a Nim sequence $G$ over $(\Y_x )$ with additive period length $\po$ that divides $Rp$, and preperiod length $\Pt$. For all $x \geq \Pt$, where $i,j \in \{ 0, \ldots, p-1 \}$, $r, r' \in \{ 0, \ldots, R-1\}$, $C, C' \in \N$ is uniquely defined by $x = j + rp + CRp$, $x-1 = i+ r'p + C'Rp$, we have $G(x)$ defined by

$$G(x) = \left\{ \begin{array}{ll}
x + d_r(j) & \mathrm{if} \: j \in S_{j,r} \\
x + d_r(i) & \mathrm{if} \:j \in T_{i,r'} \\
x & \mathrm{if} \: j \notin S_{j,r} \cup T_{i,r'} .
\end{array} \right. $$
\end{lemma}
\begin{proof}
First, $G$ is uniquely defined: If $j \in T_{i, r'}$, either i) or iii) from Definition \ref{d:DirSuc} is true, and in both cases $j \notin S_{j,r}$.

Second, $G$ is injective. If $j \notin S_{j,r} \cup T_{i,r'}$, then $j \in T_{i,r'}^*$ leads to a contradiction with iii) or ivA) in Definition \ref{d:DirSuc}, which gives $j \in T_{i,r'}$ or $j \in S_{j,r}$. Similarly, $j \in S_{j,r}^*$ leads to a contradiction with i) or ii). If $j \in S_{j,r} \cup T_{i,r'}$, injectivity is guaranteed by d) in Definition \ref{d:CutSet} and rules i)-ivB) in Definition \ref{d:DirSuc}.

Set $\Pt = \left| \min_{r,j} d_r(j) \right|$. For $x \geq \Pt$, we define the difference function $d(x) = G(x) - x$, and the cuts $( S_x, T_x)$ as in Definition \ref{d:Cut}. Then $S_x \sim S_{j,r}$, $T_x \sim T_{j,r}$, so $\# S_x = \# T_x$ for $x \geq 2 \Pt $, when we define the seed $\left[ g_0 , \ldots , g_{\Pt - 1} \right]$ as

\begin{eqnarray*}
& & \left\{ 0, \ldots , 2 \Pt - 1 \right\} \setminus G \left( \left\{ \Pt, \ldots , 3 \Pt - 1 \right\} \right) \\
& = & \left\{ 0, \ldots , 2 \Pt - 1 \right\} \setminus \left( S_{\Pt - 1}^* \cup S_{2 \Pt - 1 }^* \cup \left( G \left( \left\{ \Pt, \ldots , 2 \Pt - 1 \right\} \right) \setminus S_{\Pt - 1}^* \setminus T_{2 \Pt - 1}^*  \right) \right) ,
\end{eqnarray*}
where $S_{\Pt - 1}^*$ and $S_{2 \Pt - 1 }^*$ are disjoint, as $\Pt$ is the maximum of $d_r (j) $.

Now the additive period begins directly after the seed, as

\begin{eqnarray*}
\Pt & = & 2 \Pt - \left( \# S_{\Pt - 1} + \# S_{2 \Pt - 1 } + \left( \Pt - \# S_{\Pt - 1} - \# S_{2 \Pt - 1 } \right) \right) \\
& = & \# \left( \left\{ 0, \ldots , 2 \Pt - 1 \right\} \setminus G \left( \left\{ \Pt, \ldots , 3 \Pt - 1 \right\} \right) \right).
\end{eqnarray*}
The seed is then equal to $\N \setminus \left\{ G(x) | x \geq \Pt \right\}$. When we set $G(x) = g_x$ for $x < \Pt$, then $G$ will be injective.

Finally, $G$ is surjective. If we assume $x \in \N \setminus G( \N )$ with $G(x) > x$, we have $\# S_{j,r} = \# S_{i,r'} + 1$ and $\# T_{j,r} = \# T_{i,r'}$, which contradicts c) in Definition \ref{d:CutSet}. We reach a similar contradiction with $G(x) < x$ and $G(x) = x$.

Now we define $\Y_x$ for $x \geq \Pt$ as

\begin{eqnarray*}
\Y_x & = & \Big( \left\{ z \in \No \fw \textrm{$z$ is positively excluded at index $j$ in $\C_{j,r} $} \right\} \\
& & \cup \left\{ z \in - \No \fw \textrm{$z$ is negatively excluded at index $j$ in $\C_{i,r'} $} \right\} \\
& & \cup \left\{ z = 0 \fw j \in S_{j,r} \cap S_{j,r}^* \right\} \Big) + x ,
\end{eqnarray*}
and show that $G(x) = \mex \left( \left\{ G(x') | x' < x \right\} \cup \Y_x \right)$. Assume that $(x,y)$ is an inversion of $G$. If $G(y) - x > 0 $, then $G(y) - x $ is positively excluded in $\C_{j,r}$, so $G(y) = (G(y)-x ) + x \in \Y_x$. We find negative and zero exclusions in the same way.

It follows from b) in Definition \ref{d:CutSet2} and b) in Definition \ref{d:DirSuc2} that $d(x)$ cannot be excluded at the index $x$, so $G(x) = x + d(x) \notin \Y_x$.
\end{proof}

The most interesting conclusion from these two lemmas is that given a Nim sequence $G$ over $\left( Y_x \right)$, we can construct a seed of length no larger than $|\Mu |$ as in the previous proof. This seed will define, together with $\left( Y_x \right)$, another Nim sequence $G'$ that will have the same difference period as $G$. Of course, we could also have proven this fact directly.

\section{Optimization}%%%%%%%%%%%%%%%%%%%%%%%%%%%%%%%%%%%%%%%%%%%%%%%%%%%%%%%%%%%

In this section, we will use $\C_0$ as a standard representation of the equivalence class $[ \C_0 ]$.

While our new framework has removed the sets $(\Y_x)$, we are now forced to deal with exclusions. However, all excluded values are strictly bounded between $\Mu$ and $\Mo$, so if all difference values in a cut set are either $\Mu$ or $\Mo$, they could never be excluded. This motivates our next definition.

\begin{defn}[\textbf{Optimized cut set}]
Let $\Ch_0 = \left( \Sh_{0,r} , \Th_{0,r}, \dhat_r \right)_{r=0}^{R-1} \in \Vt \left( R, \Mu, \Mo \right)$. We call $\Ch_0$ for a \emph{optimized cut set}, if and only if $\bigcup_{r=0}^{R-1} \dhat_r \left( \Sh_r \cup \Th_r \right) \subseteq \left\{ \Mu , \Mo \right\}$. If $\C_0 = \left( S_{0,r} , T_{0,r} , d_r \right)_{r=0}^{R-1} \in \Vt$, we call $\Ch_0$ for the \emph{optimized cut set of} $\C_0$, if and only if $T_{0,r} \sim \Th_{0,r}$ and $T_{0,r}^* \sim \Th_{0,r}^*$ for all $r \in \{ 0, \ldots, R-1 \}$.

Let $\Ch_1 = \left( \Sh_{1,r} , \Th_{1,r}, \dhat_r \right)_{r=0}^{R-1} \in \Vt$ be an optimized cut set with $( \Ch_0, \Ch_1 ) \in \Et$. We call $\Ch_1$ for the \emph{optimized successor of} $\Ch_0$, if and only if $( \Ch_0, \Ch_1 )$ has no matching zeros.

In general, we call a row of a cut set $\C_{0,r}$ for \emph{optimized}, if it is optimized if considered as a cut set with one row. A difference value $d_r (x) $ is called \emph{optimized}, if $d_r (x) \in \left\{ \Mu, \Mo \right\}$.
\end{defn}

It is easy to show that the optimized cut set of $\C_0$ and the optimized successor of $\Ch_0$ are both uniquely defined.

\begin{lemma}[\textbf{Optimization Lemma}] \label{l:20}
Let $\C_0 = \left( S_{0,r} , T_{0,r} , d_r \right)_{r=0}^{R-1} \in \Vt \left( R, \Mu, \Mo \right)$. Then its optimized cut set $\Ch_0$ is a successor of $\C_0$.
\end{lemma}
\begin{proof}
We create a path $\C_0 , \ldots , \C_{M-1} , \C_M \sim \Ch_0$, where the difference values not defined by $\C_0$ are set inductively for $x = 1, \ldots, M$ by

$$ d_r ( x ) = \left\{ \begin{array}{ll}
\Mu & \textrm{if} \: x + \Mu >0 \: \textrm{and} \: x + \Mu \notin T_{x + \Mu - 1}^* \\
\Mo & \textrm{otherwise.}
\end{array} \right.$$

We show that $\C_M \sim \Ch_0$. If $x \in T_{0,r}^*$, then $d_r ( x + |\Mu | ) = \Mo$, so $M + x = |\Mu | + \Mo + x \in T_{M, r}^*$. If instead $x \notin T_{0,r}^*$, either $d_r ( x + |\Mu|) = \Mu \neq \Mo$, which implies $x + |\Mu| + \Mo =  M + x \notin T_{M, r}^*$. Or $d_r( x - \Mo ) = \Mo$, which can only happen when $\Mo < x$, which implies $M + x \notin T_{M, r}^*$.

If $x \in T_{0,r}$, then $x \leq |\Mu|$ and $d_r(x) \neq \Mo$. Then $d_r( M + x) = \Mu$, so $M + x + \Mu \leq M$, implying $M + x \in T_{M,r}$. If instead $x \notin T_{0,r}$, either $d_r( x ) = \Mo$, so $d_r (M + x ) \neq \Mu$, or $d_r( x ) = \Mu$ with $|\Mu| < x$. In both cases, we have $M + x \notin T_{M,r}$.

$\C_M$ will be optimized, as no difference value in $\C_0$ can be a part of $\C_M$.
\end{proof}
The process of creating the path $\C_0 , \ldots , \Ch_0$ is called the \emph{optimization of} $\C_0$. We can optimize an already optimized cut set this way, which implies that if there exists a connected component in $\gta$ containing only one element, it must be a non-optimized cut set.

\begin{example} \label{e:CutSet2}
Here is the optimization of the cut set from Example \ref{e:CutSet}.

$$\scaletable{ \begin{array}{ccc|cc}
+3 & \eentry & \eentry & -1 & \eentry \\
\eentry & \eentry & +2 &  -2 & \eentry \\
\eentry & \eentry & \eentry & \eentry & \eentry
\end{array}  \quad : \quad
\begin{array}{ccc|ccccc|cc}
+3 & \eentry & \eentry & -1 & +3 & +3 & -2 & -2 & -2 & \eentry \\
\eentry & \eentry & +2 &  -2 & +3 & -2 & +3 & -2 & -2 & \eentry \\
\eentry & \eentry & \eentry & +3 & +3 & -2 & -2 & -2 & \eentry & \eentry 
\end{array} }
$$
\end{example}

\begin{lemma} \label{l:20a}
Suppose $( \C_0, \C_1 ) \in \Et ( R, \Mu, \Mo )$, and let $\Ch_0 = \left( \Sh_{0,r} , \Th_{0,r}, \dt_r \right)_{r=0}^{R-1}$ be the optimized cut set of $\C_0$. Then $\Ch_0$ has a direct successor $\Ct_1 = ( \St_{1,r}, \Tt_{1,r}, \dt_r )_{r=0}^{R-1}$ with $T_{1,r} \sim \Tt_{1,r}$ and $T_{1,r}^* \sim \Tt_{1,r}^*$ for all $r \in \{ 0, \ldots, R-1 \}$ , defined by

\begin{eqnarray*}
\St_{1,r} & = & \left( \Sh_{0,r} \setminus \{ -\Mo + 1 \} \right) \cup \{ 1 \fw 1 \in S_{1,r} \} ,\\
\Tt_{1,r} & = & \left( \Th_{0,r} \setminus \{ 1 \} \right) \cup \{ k +1 \fw k \in \N , d_r( k+1 ) = -k \} ,
\end{eqnarray*}
and the possible adjoined difference values will be $\dt_r ( 1 ) = d_r (1)$, and $\dt_r (k+1) = -k$.
\end{lemma}
\begin{proof}
As $\Th_{0,r} \sim T_{0,r}$, $\Th_{0,r}^* \sim T_{0,r}^*$, and the two possible non-optimized difference values in $\Ct_1$ are the ones adjoined to $\C_1$ when it is defined as a successor of $\C_0$, it is clear that $T_{1,r} \sim \Tt_{1,r}$ and $T_{1,r}^* \sim \Tt_{1,r}^*$.

The adjoined values stem from $\C_1$ and cannot exclude themselves. An adjoined positive value must occur at index $1$ and cannot be positively excluded by the values of $\Sh_0$, and an adjoined negative value at index $k+1$ has $k+1 + d_r(k+1) = 1$ and cannot be negatively excluded by the values of $\Th_0$. Thus, $\Ct_1$ is a valid cutset.
\end{proof}

\begin{lemma} \label{l:20b}
Suppose $\C_0 \in \Vt ( R, \Mu, \Mo )$ lies in a connected component with more than one element. Then, its optimized cut set $\Ch_0$ is an ancestor of $\C_0$.
\end{lemma}
\begin{proof}
We find a path $\C_0, \ldots, C_{j-1}, C_j \sim \C_0$ for some $j > |\Mu|$. Then, we create a path $\Ch_0, \Ct_1, \ldots, \Ct_{j-1}, \Ct_j \sim \C_j$ by adjoining difference values from $\C_1, \ldots, \C_j$ to $\Ch_0$ as in Lemma \ref{l:20a}. 
\end{proof}

It follows from Lemma \ref{l:20} and \ref{l:20b} that each connected component in $\gta$ containing more that one element is fully characterized by its optimized cut sets, and the cut sets that connect them. The following lemma shows that these optimized cut sets can be connected with paths using near-optimized cut sets. The enumerations of the exceptions in the lemma should be compared with the enumerations in Definition \ref{d:DirSuc}.

\begin{lemma} \label{l:22}
Let $\left( \C_0, \C_1 \right) \in \Vt \left( R, \Mu, \Mo \right)$ with difference function $d$. Then there exists a path $\D_0 , \ldots , \D_{NM + 1}$ in $\gta$ for some $N \in \N$ with difference function $\dt$, where $\D_0 \sim \Ch_0$ and $\D_{NM + 1} \sim \Ch_1$. For all $j \in \{ 0, \ldots , NM \}$, it holds for $(\D_j , \D_{j+1} )$ that the pair has no matching zero, and for all $r \in \{ 0, \ldots, R-1 \}$, $x \in S_j \cup T_j \cup S_{j+1} \cup T_{j+1}$, we have $\dt_r (x) \in \left\{ \Mu, \Mo \right\}$ --- except that exactly one of the following exceptions apply:
\begin{itemize}
\item[(i)] There are no exceptions, so $\D_{j+1}$ is the optimized successor of $\D_j$;
\item[(ii)] There exists $x \in S_j \cup S_{j+1}$ and a difference value $k \in \{ 1, \ldots, \Mo - 1 \}$ such that for all $r \in \{ 0, \ldots , R-1 \}$, we have $\dt_r (x) \in \left\{ \Mu , k, \Mo \right\}$;
\item[(iii)] There exists $x \in T_j \cup T_{j+1}$ and a difference value $-k$ with $k \in \{ 1, \ldots, |\Mu + 1| \}$ such that for all $r \in \{ 0, \ldots , R-1 \}$, we have $\dt_r (x) \in \left\{ \Mu , -k, \Mo \right\}$;
\item[(iv)] There are matching zeros in $(\D_j , \D_{j+1} )$.
\end{itemize}
\end{lemma}
\begin{proof}
Set $N_+ = \# \left( \left\{ d_r (1) \fw r=0, \ldots, r-1 , d_r( 1 ) > 0 \right\} \setminus \left\{ \Mo \right\} \right)$, and let these positive difference values $d_r (1)$ smaller than $\Mo$ be $k_1 > k_2 > \ldots > k_{N_+ }$. \\
Set $N_- = \# \left( \left\{ d_r ( k+1 ) \fw r=0, \ldots, r-1 , k \in \No, k + d_r( k+1 ) = 0 \right\} \setminus \left\{ \Mu \right\} \right)$, and let these negative difference values larger that $\Mu$ be $-k_1' > -k_2' > \ldots > -k_{N_-}'$. Let $N_0 \in \{ 0, 1 \}$ be the indicator for whether $\left( \C_0, \C_1 \right)$ has matching zeros, and set $N = N_0 + N_+ + N_-$.

Initially, we define a path $\D_0 , \ldots ,\D_{M+1}$, where $\D_0 \sim \Ch_0$, and the $T$-- and $T^*$-- sets of $\D_1$ are equivalent to those of $\C_1$ as in Lemma \ref{l:20a}. Then we define $\D_1 , \ldots, \D_{M+1}$ as the optimization of $\D_1$. Now we extend this path between its endpoints, so it becomes

$$ \D_0 , \ldots, \D_{N_0 M + 1}, \ldots , \D_{(N_0 + N_+ )M + 1} , \ldots, \D_{NM + 1} \sim \Ch_1 .$$

All difference values in this extended path will be optimized, except for those values and matching zeros corresponding to the non-optimized values of $d$. In $\D_1$, we adjoin the matching zeroes. Between $\D_{N_0 M + 1}$ and $\D_{(N_0 + N_+ )M }$, we adjoin the positive non-optimized $d$-values. In the final part of the path, we adjoin the negative non-optimized $d$-values. The process is as follows:

\emph{(ivB)} If $N_0 = 1$, so there exist one or more matching zeros, we create $\D_{1,r}$ by adding nothing to the row, if there is a matching zero in $\left( \C_{0,r} , \C_{1,r} \right)$; if $\left( \C_{0, \tilde{r} } , \C_{1, \tilde{r} } \right)$ has no matching zeroes, we define $\D_{1,\tilde{r} }$ as the optimized successor of $\D_{0, \tilde{r} }$.

To check that this process is valid, we look at a path $\C_0', \ldots , \C_{3M + 1}'$ of cut sets with one row with difference function $d'$. First we set $\C_0' \sim \D_{0,r}$, and $\C_0' , \ldots, \C_M'$ is the optimization of $\C_0'$, so again $\C_M' \sim \D_{0,r}$. We assume there is a matching zero in row $r$ and adjoin it in $\C_{M+1}'$, and the rest of the path is repeated optimization, so $\C_{2M+1}' \sim \C_{3M + 1}' \sim \D_{ N_0 M + 1, r}$. Due to the matching zeros, we find that some difference values in this row must have a certain value:

$$\scaletable{ \begin{array}{r||c|cc|c|cc|c|c}
\textrm{Certain $d'$-value} &  \Mo & & \Mu & 0 & & \Mo & \Mu & \Mo \\
\textrm{Index} & 1 & \Mo + 1 & |\Mu | + 1 & M + 1 & M + \Mo + 1 & M + |\Mu | + 1 & 2M + 1 & 2M + |\Mu | + 1
\end{array} }$$
(Indices with blank $d'$-values cannot be determined.) For instance, as $d'( M + 1 ) = 0$, we must have $d'( |\Mu| + 1 ) \neq \Mo$ and thus $d'( |\Mu| + 1 ) = \Mu$, as the path is optimized before the index $M+1$; and as $d'( M+1 ) \neq \Mu$, so $M+1 + d'(M+1) \neq \Mo + 1$, we must have $d'( 1 ) = \Mo$, so $1 + d'( 1 ) = \Mo + 1$.

When we compare the optimized difference values in $\C_1' , \ldots , \C_M'$ with $\C_{2M + 1}' , \ldots , \C_{3M}'$, we see that for $x \in \{1, \ldots , M \}$,

$$\begin{array}{rcll}
d'( x ) & \neq & d'( 2M + x ) & \textrm{for $x = 1$, $x = |\Mu | + 1$} \\
d'( x ) & = & d'( 2M + x ) & \textrm{otherwise.}
\end{array}$$
We now examine the possible exclusions in the beginning and the end of the path. In $\C_1'$ it might happen that $k_n \in \left\{ 1, \ldots, \Mo - 1 \right\}$ is excluded at the index $1$, and with $k_n' \in \left\{ 1, \ldots, |\Mu + 1| \right\}$, we might exclude $-k_n'$ at the index $k+1$. This is why we need to change the rows with the matching zeros first in the path $\D_0, \ldots, \D_{NM + 1}$. After the optimization of $\C_{M+1}'$, nothing can be excluded at the index $2M + 1$. There might appear new exclusions at $M + |\Mu | + 1$, but as no new non-optimized difference values are adjoined at the index $aM + |\Mu | +1$ for $a \in \N$ in $\D_{N_0 + 1}, \ldots, \D_{NM + 1}$, these exclusions do not matter.

\emph{(ii)} For $n \in \{ 1, \ldots, N_+ \}$, we create $\D_{(N_0 + n - 1)M + 1,r}$ by adjoining $k_n$ to the row $r$ at the index $(N_0 + n - 1 )M + 1$, if the corresponding row in $\C_1$ fulfills $d_r( 1 ) = k_n$; otherwise, we define $\D_{(N_0 + n - 1)M + 1,\tilde{r} }$ as the optimized successor of $\D_{(N_0 + n - 1 )M,\tilde{r} }$. We reuse the notation $\C_0', \ldots , \C_{3M + 1}'$, where we now have $\C_0' \sim \D_{(N_0 + n - 1 )M,r}$, and $\C_0' , \ldots, \C_M'$ is the optimization. We adjoin the difference value $k_n$ in $\C_{M+1}'$, and the rest of the path will be optimized. Now, the certain difference values become:

$$\scaletable{ \begin{array}{r||c|cc|c|cc|c|c}
\textrm{$d'$-value} & \Mo & & \Mu & +k_n & & \Mo & \Mu & \Mo \\
\textrm{Index} & 1 &\Mo \! + \! 1 & |\Mu | \! + \! k_n  \! + \! 1 & M \! + \! 1 &  M \! + \! \Mo \! + \! 1 & M \! + \! |\Mu | \! + \! k_n \! +  \! 1 & 2M \! + \! 1 & 2M \! + \! |\Mu | \! + \! k_n \! + \! 1
\end{array} }$$
For $x \in \{1, \ldots , M \}$,

$$\begin{array}{rcll}
d'( x ) & \neq & d'( 2M + x ) & \textrm{for $x = 1$, $x = |\Mu | + k_n + 1$ } \\
d'( x ) & = & d'( 2M + x ) & \textrm{otherwise.}
\end{array}$$
It might happen that $k_{n'}$ is excluded at the index $1$, if $k_{n'} > k_n$. If $k_{n'} < k_n$, then $k_{n'}$ cannot be exluded at the index $1$ or $M+1$, as then $\C_{1,r}$ could not be a successor of $\C_{0,r}$. The possible new exclusions at $M + |\Mu | + k_n +  1 $ cannot affect the adjoined values in $\D_0 , \ldots, \D_{NM + 1}$.	

\emph{(iii)} For $n' \in \{ 1, \ldots, N_- \}$, we create $\D_{(N_0 + N_+ + n' - 1)M + 1,r}$ by adjoining $-k_n'$ to the row $r$ at the index $(N_0 + N_+ + n' - 1 )M + k_n' + 1$, if the corresponding row in $\C_1$ has the difference value $d_r( 1 + k_n') = -k_n'$; otherwise, the row becomes the optimized successor as above. We redefine the path $\C_0', \ldots , \C_{3M + 1}'$ similarly, so the certain difference values become:

$$\scaletable{ \begin{array}{r||c|cc|c|cc|c|cc}
\textrm{$d'$-value} & & \Mo & \Mu &  & -k_n' & \Mo  & & \Mu & \Mo \\
\textrm{Index} & 1 &k_n' \! + \! 1 & |\Mu | \! + \! 1 & M \! + \! 1 &  M \! + \! k_n' \! + \! 1 & M \! + \! |\Mu | \! +  \! 1 & 2M \! + \! 1 & 2M \! + \! k_n' \! + \! 1 & 2M \! + \! |\Mu | \! + \! 1 
\end{array} }$$
For $x \in \{1, \ldots , M \}$,

$$\begin{array}{rcll}
d'( x ) & \neq & d'( 2M + x ) & \textrm{for $x = k_n' + 1$, $x = |\Mu | + 1$ } \\
d'( x ) & = & d'( 2M + x ) & \textrm{otherwise.}
\end{array}$$
It might happen that $-k_{n'}'$ is excluded at the index $1 + k_n'$, if $k_{n'}' > k_n'$. The possible new exclusions at $M + |\Mu | +  1 $ cannot affect the adjoined values in $\D_0 , \ldots, \D_{NM + 1}$.

\emph{(ivA)} We should mention the case when two difference values $k$ and $-k'$ are adjoined to $\C_{1,r}$. If either of these difference values are optimized, we proceed as above. If both are non-optimized, so $k = k_n$ and $k' = k_{n'}'$, we adjoin $k_n$ to $\D_{(N_0 + n - 1)M + 1,r}$, and $-k_{n'}'$ to $\D_{(N_0 + N_+ + n' - 1)M + 1,r}$; all other values of $\dt_r$ are optimized. The fixed difference values will then be as in case (ii) and (iii) combined.

For all $j \in \{ 0, \ldots, (N + 1)M \}$ and $r, r' \in \{ 0, \ldots, R-1 \}$ with $r \neq r'$, we must have $\D_{j,r} \not \sim \D_{j, r'}$. If we assume that $\D_j$ is an invalid cut set with $\D_{j,r} \sim \D_{j, r'}$, Lemma \ref{l:13c} shows that for all $j'>j$, $\D_{j'}$ would be invalid for the same reason, but $\D_{NM+1}$ is a valid cut set.
\end{proof}
In passing, we have shown that the order of the matching zeros and the non-optimized difference values in the path cannot be altered.

\begin{example} \label{e:CutSet3}
We create a near-optimized path in $\gta (3, -2, 3 )$ between the cut sets from Example \ref{e:CutSet}.

$$\scaletable{ \begin{array}{ccc|cccccccccc|cc}
+3 & \eentry & \eentry & -2 & +3 & +3 & -2 & -2 & -2 & +3 & +3 & -2 & -2 & -2 & \eentry \\
\eentry & +3 & \eentry & -2 & +3 & -2 & +3 & -2 & -2 & -1 & +3 & +3 & -2 & -2 & -2 \\
\eentry & \eentry & \eentry & +2 & +3 & -2 & -2 & +3 & -2 & -1 & +3 & -2 & +3 & -2 & -2 
\end{array} }$$
\end{example}

\begin{remark} \label{rem:22}
We notice as a general pattern in the $\C'$--paths that before a non-optimized difference value is adjoined, it is excluded at the same index $\pmod M$. In case (ii), if a positive difference value $k_n$ is adjoined to the near-optimized path at index $(N_0 + n - 1)M + 1$, it will be excluded at the indices $aM + 1$ with $a < N_0 + n - 1$. Similarly, before $-k_{n'}'$ is adjoined , it will be excluded at the indices $aM + 1 + k_{n'}'$ for $a < N_0 + N_+ + n' - 1$. A matching zero is excluded at index $1$ before it is adjoined.

It follows that if $\Ch_0 , \ldots , \Ch_M $ is a path of optimized cut sets, and $k_n$ is adjoined to $\C_{M+1, r}$, then in all rows $r'$ where $k_n$ is excluded at the index $1$, the value must be adjoined to $\C_{M+1, r'}$. The same is true for matching zeros. For negative difference values, if $-k_{n'}'$ is adjoined to $\C_{M+1, r}$ at the index $M + 1 + k_{n'}'$, the value must be adjoined in all other rows where $-k_{n'}'$ was excluded at the index $1 + k_{n'}'$.
\end{remark}

\section{Binary representations} %%%%%%%%%%%%%%%%%%%%%%%%%%%%%%%%%%%%%%%%%%%%%%%%%%%%%%%%%%

In this section we only need to deal with optimized difference values, so we simplify the notation by writing $\Pl$ for $\Mo$, and $\Mi$ for $\Mu$.

\begin{defn}[\textbf{Binary representation}]
Let $\C_0 \in \Vt \left( R, \Mu, \Mo \right)$. The \emph{binary re\-presentation of} $\C_0$ is an $R \times M$ matrix $\B_0 = \left[ b_{0,r,m} \right]_{r=0, \ldots, R-1 , m=0, \ldots, M-1}$ with entries in $\{ \Pl, \Mi \}$ defined by

$$b_{0, r, m} = \left\{ \begin{array}{rl}
\Pl & \textrm{if $m - \Mu + 1 \in \Sh_{0,r}$ for $m < \Mo$, or $m - \Mu + 1 \notin \Th_{0,r}$ for $m \geq \Mo$ } \\
\Mi & \textrm{if $m - \Mu + 1 \in \Th_{0,r}$ for $m \geq \Mo$, or $m - \Mu + 1 \notin \Sh_{0,r}$ for $m < \Mo$ ,} 
\end{array} \right. $$
where $\Ch_0 = \left( \Sh_{0,r} , \Th_{0,r}, \dhat_r \right)_{r=0}^{R-1}$ is the optimized cut set of $\C_0$.

Let $c_m$ be the symbol for the $m$'th column of $\B_0$, so

$$ c_m = \left[ \begin{array}{c} 
b_{0,0,m} \\ \vdots \\ b_{0, R-1 ,m }
\end{array} \right] .$$

Let $\tau $ be the column permutation that cycles all columns of $\B_0$, defined by $\tau (m) = ( m - 1 ) \bmod M$.
\end{defn}

The permutation $\tau$ represents the rotation of the columns from $\B_0$ to $\B_1$, where $\B_1$ is the binary represantion of a direct successor $\C_1$ of $\C_0$. For instance, if $\Ch_1$ is the optimized successor of $\Ch_0$, then $\B_{0, r, m } = \B_{1, r, \tau (m) } $ for all $r, m$.

Note that $\B_0$ cannot have two identical rows, as it would violate condition a) in Definition \ref{d:CutSet2}.
\begin{defn}[\textbf{Reversible}]
Let $\left( \Ch_0 , \C_1 \right) \in \Et \left( R, \Mu, \Mo \right)$ fulfill the conditions from Lemma \ref{l:22}. We call $\left( \Ch_0 , \C_1 \right) $ for \emph{reversible}, if and only if $\Ch_0$ and $\C_1$ lie in the same connected component. If this holds, we say that there is a \emph{reversible shift} in their binary representations from $\B_0$ to $\B_1$.
\end{defn}

\begin{example}
Continuing Example \ref{e:CutSet}, here are the binary representations of the cut set and its successor.

$$\scaletable{ \begin{array}{ccc|cc}
\Pl & \Mi & \Mi & \Mi & \Pl \\
\Mi &\Pl  & \Mi & \Mi & \Pl  \\
\Mi & \Mi & \Mi & \Pl  & \Pl 
\end{array} \quad \hookrightarrow \quad \begin{array}{ccc|cc}
\Mi & \Mi & \Mi & \Pl  & \Pl \\
\Pl  & \Mi & \Mi & \Mi & \Pl \\
\Mi  & \Pl & \Mi & \Mi & \Pl  
\end{array} }$$
When we ignore the $\tau$--rotation of the columns, we see that $\left( c_0 , c_4 \right)$ and $\left( c_2 , c_3 \right)$ in the ancestor have been swapped. As we will prove below, all these shifts are reversible.
\end{example}

\begin{lemma} \label{l:23}
Let $\left( \Ch_0 , \C_1 \right) \in \Et \left( R, \Mu, \Mo \right)$ fulfill the conditions from Lemma \ref{l:22}. Then the following two statements are equivalent:
\begin{itemize}
\item[a)] $\left( \Ch_0 , \C_1 \right)$ is reversible.
\item[b)] Either $\B_0 / \tau = \B_1 / \tau$, or we can create $\B_1 / \tau$ by swapping two columns in $\B_0 / \tau$.
\end{itemize}
\end{lemma}
\begin{proof}
\emph{(i)} If $\C_1$ is the optimized successor of $\Ch_0$, then $\B_0 / \tau = \B_1 / \tau$. To create a path from $\C_1$ to $\Ch_0$, we optimize $\C_1$.

\emph{(ii), (iii), (ivB)} In the proof of Lemma \ref{l:22}, we showed that in each of these cases, exactly two of the optimized difference values at the indices modulo $M$ will change. So for each row $r$ where a non-optimized difference value is adjoined, a $\Pl$ is swapped with a $\Mi$:

$$\begin{array}{ll|ccccc|ccccc|}
& & \multicolumn{4}{c}{\textrm{In $\B_0 / \tau$ at column: }} & & \multicolumn{4}{c}{\textrm{In $\B_1 / \tau$ at column:}} & \\
\textrm{Case} & \textrm{Value adjoined} & 0 & k & |\Mu| & k + |\Mu| & & 0 & k & |\Mu| & k + |\Mu| & \\
\hline
(ii) & +k & \Pl &     &       & \Mi & & \Mi &     &     & \Pl & \\
(iii) & -k  &     & \Pl & \Mi &      & &       & \Mi & \Pl &    & \\
(ivB) & [0] & \Pl &  & \Mi &  &  & \Mi &  &  \Pl &  & \\
\end{array}$$
(We have normalized the column indices, so $+k$ will be adjoined at the index $0$.)

It follows that if $b_{0,r,m_1} = \Pl$, $b_{0,r,m_2} = \Mi$ with $m_1, m_2 \in \{ 0 , \ldots, M-1 \}$, we can always find a successor $\C_j$ to $\Ch_0$ where these two entries become swapped in $\B_j / \tau$. Begin to optimize $\Ch_0$, then adjoin to the path the following non-optimized difference value in column $c_{m_1}$:

$$\textrm{If  $( m_1 - m_2 ) \bmod M$ } \left\{ \begin{array}{ll}
= |\Mu| & \textrm{adjoin a matching zero.} \\
< |\Mu| & \textrm{adjoin $( (m_2 - m_1 - |\Mu| ) \bmod M ) - M $. }\\
> |\Mu| & \textrm{adjoin $(m_2 - m_1 - |\Mu| ) \bmod M$. } \\
\end{array} \right.$$
We can do the same with $\Ch_1$, so when we swap $b_{0,r, m_1 }$ and $b_{0,r, m_2}$, we can create a path that will swap these entries back. The question is what happens in the other rows.

Let $r', r'' \in \{ 0, \ldots, R-1 \}$. By Remark \ref{rem:22}, if $b_{0, r', m_1} = \Pl $ and $b_{0, r', m_2} = \Mi$, these entries will be swapped as well. Suppose that $\left( b_{0, r'', m_1 }, b_{0, r'', m_2 } \right) \neq \left( \Mi, \Pl \right)$ for all $r''$. Then we are swapping the entire columns $c_{m_1}$ and $c_{m_2}$, and when we swap the entries back in the $r$'th row, all other entries will be swapped. So $\left( \Ch_0 , \C_1 \right)$ is reversible.

However, if $\left( b_{0, r'', m_1 }, b_{0, r'', m_2 } \right) = \left( \Pl , \Mi \right)$ for some $r''$, the swap will result in $\left( c_{m_1}, c_{m_2} \right)$ being changed, as they have different numbers of pluses and minuses after the swap. So when the entries of the $r$'th row are swapped back to their original state, $b_{0, r'', m_2}$ and $b_{0, r'', m_1}$ will also be swapped, so $\left( c_{m_1}, c_{m_2} \right)$ will remain changed after the two shifts. It does not matter which path we follow and how many other columns we shift in the binary representation, as the entries $b_{0, r, m_1}$ and $b_{0, r'', m_1}$ now shift together. Thus, $\left( \Ch_0 , \C_1 \right)$ becomes irreversible.
\end{proof}
It follows that if $\C_x$ and $\C_y$ lie in the same connected component, then we can create $\B_y$ by some column permutation of $\B_x$. Using this fact, we can prove our main theorem.

\begin{thm}\label{t:Main3}
Any Nim sequence $G$ is additively periodic, and the length $\po$ of its period is bounded by $\po \leq K_{\Mu , \Mo}\:  p$. Here, $K_{\Mu , \Mo}$ is the maximal value of $\lcm \left( p_1, \ldots, p_n \right)$, where $p_1, \ldots , p_n \in \No$ for some $n \in \No$ with the constraints that $\sum_{i=1}^n p_i = M - 1$, and $n \leq \min \left( |\Mu| , \Mo \right)$.

Also, there exists a Nim sequence $G$ for which $\po = K_{\Mu , \Mo} \: p$.
\end{thm}

\begin{proof}
Given a Nim sequence $G$ with additive period length $\po$, we set $R = \frac{ \lcm \left( \po, p \right) }{p}$ and create the path $\C_0, \ldots, \C_p$ as in Lemma \ref{l:17}, where $\C_p$ is the cycled cut set of $\C_0$. As we can continue the path to reach $\C_{Rp} \sim \C_0$, it follows that $\B_p$ is created by a column permutation of $\B_0$. Let $c_0 , \ldots, c_{M-1}$ be the columns of $\B_0$, and let $\pi \left( c_m \right) = \left[ b_{0, \pi (r), m}  \right]_{r=0}^{R-1}$. It follows from Definition \ref{d:CycCut} that $\pi \left( c_0 \right), \ldots, \pi \left( c_{M-1} \right)$ is a permutation of $c_0 , \ldots, c_{M-1}$.

We say that $\left\{ c_m, \pi \left( c_m \right), \ldots, \pi^{p_i - 1} \left( c_m \right) \right\}$ constitute a $p_i$--cycle, if $p_i \in \No$ is the smallest number such that $\pi^{p_i} \left( c_m \right) = c_m$. The first elements $\left[ b_{0, r, m} \right]_{r=0}^{p_i - 1}$ will be repeated in $c_m$, so $b_{0, r, m} = b_{0, r + a p_i , m}$ for $a \in \{ 0, \ldots, R / p_i - 1 \}$. Now the elements of $c_m$ are all identical, if and only if $\left\{ c_m \right\}$ is a $1$--cycle. If there are $n$ column cycles of length $p_1, \ldots, p_n$ in $\B_0$, we must have $R = \lcm \left( p_1, \ldots, p_n \right)$, because every $p_i$ must divide $R$, and two rows in $\B_0$ cannot be identical.

We need to examine how the columns shift from $\B_0$ to $\B_p$. If $p_i$ happens to divide $M$, we can space out the columns of the $p_i$--cycle at the indices $0, R / p_i  , \ldots, R - R / p_i$, and use the $\tau$--rotation to shift the columns; otherwise, we need to swap at least one column of the $p_i$--cycle. If $p_i \geq 2$, and $c_{m_1}, c_{m_2}$ belong to the same $p_i$--cycle, it is easy to see that there exist row indices $r', r''$ such that $\left( b_{0, r', m_1} , b_{0, r', m_2} \right) = \left( b_{0, r'', m_2} , b_{0, r'', m_1} \right) = ( \Pl, \Mi)$, so we cannot swap these columns by Lemma \ref{l:23}. If $c_{m_1}$ belongs to a $p_i$--cycle and $c_{m_2}$ belongs to a $p_{i'}$--cycle with $p_i$, $p_{i'}$ coprime, we cannot swap $c_{m_1}$ and $c_{m_2}$ for the same reason. However, if $\left\{ c_{m_1} \right\}$ is a $1$--cycle, we can swap $c_{m_1}$ with any another column in a reversible shift.  To ensure a reversible shift is possible from $\B_0$ to $\B_p$, we just need to include a $1$--cycle in $\B_0$.

The problem of maximizing $R$ is equal to finding $\max \left( \lcm \left( p_1, \ldots, p_n \right) \right)$ under the constraint that $\sum_{i=1}^n p_i = M - 1$. Finally, there is a bound for $\Mu$ and $\Mo$. As each of the $n$ cycles must contain both one $\Pl$ and one $\Mi$, we must have $\min \left( |\Mu| , \Mo \right) \geq n$, before we can achieve the maximum value for $R$.

A Nim sequence $G$ with maximal period length given $\Mu$ and $\Mo$ can be found by first creating a binary representation which permits the desired number of $p_i$--cycles with lengths $1 , p_1 , \ldots , p_n$, which are found under the constraints given above. For each column swap, we find the resultating binary representation, convert them to optimized cut sets, and connect them with a near-optimized path. Finally, we use Lemma \ref{l:18} to find $\left( \Y_x \right)$ and the seed.
\end{proof}

\begin{proof}[Proof of Theorem \ref{t:Main}]
To maximize $R$, we need to have as many prime factors as possible in $p_1, \ldots, p_n$. So if $M-1$ is equal to the sum of the $n$ smallest prime numbers, the solution is $\prod_{i=1}^n p_i$, where $p_1 = 2$, $p_2 =3$, \ldots, and $p_n$ is the $n$'th smallest prime.

By an extension of the Prime Number Theorem (see \cite{PrimeSum}, \cite{Salat}), the asymptotic value of the sum of all primes smaller or equal to $x \in \No$ is

$$ \sum_{p_i \leq x } p_i \sim \Li \left( x^2 \right)  \quad \textrm{as $x \to \infty$ },$$
and the product of all primes smaller or equal to $x$, the \emph{primorial} (see \cite{Primorial}), has the asymptotic value

$$ \prod_{ p_i \leq x} p_i \sim \exp (x) \quad \textrm{as $x \to \infty$ } .$$
Combining these formulas, we find

$$\max \left( \lcm \left( p_1, \ldots, p_n \right) \right) \sim \exp \left( \sqrt{ \Li^{-1} \left( M-1 \right) } \right) \quad \textrm{as $M \to \infty$ } $$
which gives the asymptotic value of $K_{\Mu , \Mo }$.
\end{proof}

\begin{remark}
The proof assumes that we can choose the elements of each $\Y_x$ freely between the difference bounds. If $\left( \Y_x \right)$ is given, the maximal bound $K_{\Mu, \Mo} p$ may shrink, if the excluded elements do not permit the optimal number of column cycles. The bounds of Theorem \ref{t:Main2} and Remark \ref{rem:Main} may also shrink.
\end{remark}

\begin{example} \label{e:gtaM6}
Let us find the maximum value of $R$ in $\gta ( R, -3, 3 )$. As $M=6 = 1+2+3$, we define a binary representation with a $1$-cycle, a $2$-cycle and a $3$-cycle (example, left below). As it happens that $\max R = 2*3 = 6$, we could define a binary representation as in Example \ref{e:Simple} with one $6$-cycle, which cycles its columns by the $\tau$--rotation (example, right below). In this case, $\{ -2, \ldots, 2 \}$ are excluded at every index, so we must have $p=1$, and the difference period will be $[ +3, +3, +3, -3, -3, -3]$.

$$\scaletable{ \begin{array}{ccc|ccc}
\Pl & \Pl & \Mi & \Pl & \Mi & \Mi \\
\Pl  &\Mi  & \Pl & \Mi & \Pl & \Mi \\
\Pl  & \Pl & \Mi & \Mi & \Mi & \Pl \\
\Pl  & \Mi & \Pl & \Pl & \Mi & \Mi \\
\Pl  & \Pl & \Mi & \Mi & \Pl & \Mi \\
\Pl  & \Mi & \Pl & \Mi & \Mi & \Pl 
\end{array} \quad ; \quad \begin{array}{ccc|ccc}
 \Pl & \Pl & \Pl & \Mi & \Mi & \Mi \\
 \Mi & \Pl & \Pl & \Pl & \Mi & \Mi \\
 \Mi & \Mi & \Pl & \Pl & \Pl & \Mi \\
\Mi  & \Mi & \Mi & \Pl & \Pl & \Pl \\
\Pl  & \Mi & \Mi & \Mi & \Pl & \Pl \\
\Pl  & \Pl & \Mi & \Mi & \Mi & \Mi 
\end{array} }$$

If we look at $\gta ( R, -4, 2 )$, we can use the left binary representation if we change the $1$-cycle from all pluses to all minuses. Changing the right example, we can define another $6$-cycle so the difference period becomes $[+2, +2, +2, +2, -4, -4]$. If we look at $\gta (R, -5, 1)$, we can no longer use the left example, but we can define a $6$-cycle so the difference period becomes $[+1, +1, +1, +1, +1, -5]$.
\end{example}

We now expand our framework, so we can give a bound for the length of the preperiod, or rather $\Pt - L$, where $L$ is the length of the seed. This bound must also depend on the seed, as we can make the preperiod length arbitrarily long by including large elements in the seed.

Given $\Ro, R \in \No$, we now consider a cut set $\C_x = \left( S_{x,r} , T_{x,r} , d_r \right)_{r=0}^{\Ro + R - 1}$ to have $\Ro + R$ rows, where the top $\Ro$ rows define the preperiod. For $r \in \{ \Ro, \ldots, R + \Ro - 1 \}$, we redefine the row rotation of the period as $\pi( r ) = ( ( r - \Ro - 1 ) \bmod R ) + \Ro $, so it cycles the bottom $R$ rows. For $r \in \{ 0, \ldots, \Ro - 1 \}$, we define the row rotation of the preperiod as

$$\pio( r ) = ( r-1 ) \bmod \Ro .$$
This rotation cycles all rows with indices $0, \ldots, \Ro - 2$, as the bottom row of the preperiod will continue into the period.

A Nim sequence with preperiod length $\left( \Ro - 1 \right) p + x'$ and period length $Rp$, where $x' \in \left\{ 1, \ldots, p \right\}$, corresponds to a path $\C_0, \ldots, \C_{p-1}$ in $\gta \left( \Ro + R, \Mu, \Mo \right)$ . For $j < p - x'$, the cut sets $\C_j = \left( S_{j,r} , T_{j,r} , d_r \right)_{r=1}^{\Ro + R - 1}$ lack one top row compared to the remaining cut sets in the path. This $0$'th row will be added in $\C_{p - x'}$. One cut set $\C_p$, a direct successor of $\C_{p-1}$, is adjoined at the end of the path, where

$$\C_{0,r} \sim \left\{ \begin{array}{rl}
\C_{p, \pio( r )} & \textrm{for $r \in \{ 1, \ldots, \Ro - 1 \}$} \\
\C_{p, \pi( r )} & \textrm{for $r \in \{ \Ro, \ldots, \Ro + R - 1 \}$} \\
\C_{p, \Ro - 1} & \textrm{for $r = \Ro $} . 
\end{array} \right.$$
We have $\C_{p, \Ro - 1}  \sim \C_{p, \Ro + R - 1}$, where the preperiod shifts into the period, so $\C_p$ will violate condition a) in Definition \ref{d:CutSet2}. We can define a legal cut set $\Ct_p$ by removing the illegal row at the index $\Ro - 1$ before we move all other rows in the preperiod one index down. When now $x' < p$, $\Ct_p \sim \C_0$; or with $x' = p$, the cut sets will be equivalent when we ignore the top row of $\C_0$.

\begin{example}
The first values in the Nim sequence $G_3$ from Wythoff's Game in Example \ref{e:Wythoff} corresponds to a path $\C_0, \ldots, \C_6$ in $\gta ( 5, -5, 3)$.

$$\scaletable{ \begin{array}{ccc|cccccc|ccccc}
  &  &  &  & +3 & +3 & +3 & +3 & -2 & -5 & -5 & \eentry & \eentry & \eentry \\
\eentry & +3 & +3 & +3 & +3 & -2 & -5 & -5 & +2 & \eentry & \eentry & -4 & \eentry & \eentry \\
+3  & +3 & \eentry & -5 & -5 & +2 & +2 & +3 & -2 & -4 & \eentry & \eentry & \eentry & \eentry \\
\hline
\eentry & +3 & \eentry & +2 & +3 & -2 & -4 & +3 & -2 & \eentry & \eentry & \eentry & -4 & \eentry \\
\eentry & +3 & \eentry & -4 & +3 & -2 & +2 & +3 & -2 & -4 & \eentry & \eentry & \eentry & \eentry    
\end{array} }$$
(The horizontal line separates the preperiod from the period.)
\end{example}

\begin{thm} \label{t:24}
Let $G$ be a Nim sequence over $\left( \Y_x \right)_{x=0}^\infty$ that has period length $p$. Let $G$ have seed $\left[ g_0, \ldots, g_{L-1} \right]$, period length $\po$ and preperiod length $\Pt$.

Set $\Kh = \max_{x<L} g_x - L$. Then the preperiod length is bounded by

$$\Pt - L \leq \left( \frac{M}{2} \right)^2 p + \Kh ,$$
and for $M \geq 11$, we have

$$\Pt - L + \po \leq K_{\Mu, \Mo} \: p + \Kh,$$
where $K_{\Mu, \Mo}$ is defined as in Theorem \ref{t:Main3}.
\end{thm}
\begin{proof}
Let $c_0 , \ldots, c_{M-1}$ be the columns of $\B_0 / \tau$. We say that $\pio (c_{m_1} ) = c_{m_2}$ in $\B_j / \tau$ if and only if $b_{j, r, m_1} = b_{j, \pio (r), m_2}$ for all $r \in \{ 0, \ldots, \Ro - 1 \}$, ignoring the values of the period. Then $\left\{ c_m, \pio \left( c_m \right), \ldots, \pio^{p_i - 1} \left( c_m \right) \right\}$ constitute a $p_i$--cycle in the preperiod, if $p_i \in \No$ is the smallest number such that $\pio^{p_i} \left( c_m \right) = c_m$.

Any shift of the columns in $\B_p$ affects both the period and the preperiod. It follows that if we have a $p_i$--cycle in the preperiod, the columns of the cycle must correspond to $p_{i'}$--cycle in the period, where either $p_{i'}$ divides $p_i$, or $p_i$ divides $p_{i'}$, or $p_i = p_{i'}$.

As written above, the difference values of $G$ defines a path $\C_0, \ldots, \C_p$, where $\C_p$ has two equivalent rows. As $\C_{p-1}$ is a legal cut set, and the shift from $\B_{p-1}$ to $\B_p$ creates two identical rows in $\B_p$, it follows that this shift is irreversible. We find two columns $c_{m_1}, c_{m_2}$ where this shift swaps $b_{p-1, \Ro - 1, m_1} = \Pl$ with $b_{p-1, \Ro - 1, m_2} = \Mi$, so these entries become equal to $\left( b_{p, \Ro +R - 1, m_1} , b_{p, \Ro + R - 1, m_2} \right) = \left( \Mi, \Pl \right)$. This shift cannot affect any other rows, as the shift will also be irreversible in these rows, which should all cycle back to the entries of $\B_0$. So $\left( b_{p-1, r, m_1} , b_{p-1, r, m_2} \right) \neq \left( \Pl , \Mi \right)$ for all $r \neq \Ro - 1$. (There might be other irreversible shifts affecting the preperiod in the path, but they too cannot affect any other rows than the row $\Ro-1$.)

Assume there are two column cycles in the preperiod of length $p_1$ and $p_2$ with $p_1 , p_2$ coprime, and that all other cycles have length $p_i$ which either divides $p_1$ or $p_2$. We can define the columns of $\B_{p-1} / \tau$ in the $p_1$--cycle by the left pattern below, and define the columns in the $p_2$--cycle by the right pattern below.

$$ [ \underbrace{ \underbrace{\Mi, \Mi, \ldots, \Mi}_{\textrm{$p_1 - 1$ times} }, \Pl  }_{\textrm{repeat $p_2$ times} } ] \quad ; \quad
[ \underbrace{ \underbrace{\Pl, \Pl, \ldots, \Pl}_{\textrm{$p_2 - 1$ times} }, \Mi  }_{\textrm{repeat $p_1$ times} } ] .$$
Let $c_{m_1}$ belong to the $p_1$--cycle and $c_{m_2}$ to the $p_2$--cycle. Then there is exactly one row index $r$, where $\left( b_{p-1, r, m_1} , b_{p-1, r, m_2} \right) = \left( \Pl, \Mi \right) $, and we define the columns, so $r = \Ro - 1$. We then let $R=1$, and set $\left( b_{p-1, \Ro, m_1} , b_{p, \Ro, m_2} \right) = \left( \Mi, \Pl \right)$, and $ b_{p-1, \Ro, m} = b_{p-1, \Ro - 1, m}$ for $m \in \{ 0, \ldots, M-1 \} \setminus \left\{ m_1, m_2 \right\}$. Now the swapping of $\left( b_{p, r, m_1} , b_{p, r, m_2} \right)$ is irreversible, and makes the rows $\Ro-1$ and $\Ro$ identical.

As we must have a common $1$--cycle in the preperiod and period to swap the columns from $\B_0 / \tau$ to $\B_p / \tau$, we find the maximal value of $\Ro$ as $\Ro = p_1 p_2$ under the constraint $p_1 + p_2 = M - 1$. So if say $M$ is even, we have $\max ( p_1 p_2 ) = \frac{M}{2} \left( \frac{M}{2} - 1 \right) \to \left( \frac{M}{2}^2 \right)$ as $\Mu \to - \infty, \Mo \to \infty$.

If we have three column cycles in the preperiod of lengths $p_1, p_2, p_3$, all coprime, when $\left( b_{p-1, r, m_1} , b_{p-1, r, m_2} \right) = \left( \Pl, \Mi \right) $ for some $r < \Ro$, there will always exist $r' < \Ro$, $r' \neq r$ with $\left( b_{p-1, r', m_1} , b_{p-1, r', m_2} \right) = \left( \Pl, \Mi \right) $. Thus, it is impossible to make an irreversible shift that only affects one row in the preperiod, which shows that we cannot have more than two column cycles with coprime period lengths.

To show that $\Pt - L + \po \leq K_{\Mu, \Mo} \: p$, it is enough to examine the extreme case where $\po = \max \left( \lcm \left( p_1, \ldots, p_n \right) \right) p$ with $\sum_{i=1}^n p_i = M - 1$. If we have two columns $c_{m_1}, c_{m_2}$ that either belong to a $p_i$-- and a $p_{i'}$--cycle, or the same $p_i$--cycle, it is impossible to find an irreversible shift between $c_{m_1}$ and $c_{m_2}$ in the preperiod that would not affect rows in the period. As the $1$--cycle that is used to swap columns must have identical elements in the period and the preperiod, we cannot find an irreversible shift that affects this $1$--cycle. This assumes that we can fit three cycles of coprime period lengths in the period, and $M = 11 = 1 + 2 + 3 + 5$ is the smallest value where this is possible. So with $M \geq 11$ and $\po$ maxed out, we cannot have any values in the preperiod at all, unless we have large elements in the seed.

We can adjoin some additional $\Kh$ values to the start of the preperiod, but only if they are not excluded or cause any exclusions. So we must adjoin negative difference values equal or smaller to $\Mu$. The maximal number of negative values we can adjoin this way is $\max_{x<L} g_x - L$, as shown in Lemma \ref{l:6}.
\end{proof}

\begin{example}
We can expand the left binary representation of Example \ref{e:gtaM6} to seven rows, where the top six rows represents the preperiod. Here we have $7 = \Ro + R > 6 = K_{\Mu, \Mo}$, which is possible, as $M = 6 < 11$.

$$\scaletable{ \begin{array}{ccc|ccc}
\Pl & \Pl & \Mi & \Pl & \Mi & \Mi \\
\Pl  &\Mi  & \Pl & \Mi & \Pl & \Mi \\
\Pl  & \Pl & \Mi & \Mi & \Mi & \Pl \\
\Pl  & \Mi & \Pl & \Pl & \Mi & \Mi \\
\Pl  & \Pl & \Mi & \Mi & \Pl & \Mi \\
\Pl  & \Mi & \Pl & \Mi & \Mi & \Pl \\
\hline
\Pl  & \Pl & \Pl & \Mi & \Mi & \Mi 
\end{array} \quad \hookrightarrow \quad \begin{array}{ccc|ccc}
\Pl & \Pl & \Mi & \Pl & \Mi & \Mi \\
\Pl  &\Mi  & \Pl & \Mi & \Pl & \Mi \\
\Pl  & \Pl & \Mi & \Mi & \Mi & \Pl \\
\Pl  & \Mi & \Pl & \Pl & \Mi & \Mi \\
\Pl  & \Pl & \Mi & \Mi & \Pl & \Mi \\
\Pl  & \Pl & \Pl & \Mi & \Mi & \Mi \\
\hline
\Pl  & \Pl & \Pl & \Mi & \Mi & \Mi 
\end{array} }$$
Left is $\B_{p-1} / \tau$, right is $\B_p / \tau$. The irreversible shift swaps $b_{p-1, 5, 1}$ with $b_{p-1, 5, 5}$.
\end{example}

\begin{example}
Here is a path in $\gta ( 3, -1, 2 )$ with $p = 3$, where the preperiod is extended by $\Kh = 2$.

$$ \scaletable{ \begin{array}{cc|ccc|c}
  &  & [+2] & -1  & -1 & \eentry \\
\eentry & \eentry & +1 & -1 & +1 & -1 \\
\hline
\eentry & +1 & -1 & +2 & -1 & -1 
\end{array} }$$
As $d_0 (0) = 2$ excludes $1$, this cannot be a legal path unless $d_0 (0) $ belongs to the seed. 
\end{example}

\section{Conclusion}
%%%%%%%%%%%%%%%%%%%%%%%%%%%%%%%%%%%%%%%%%%%%%%%%%%%%%%%%%%%%%%%%

In this paper, we have established bounds for the period length and the preperiod length of Nim sequences. To do this, we defined several concepts, such as cut sets and binary representations. These concepts could possibly be used in the study of specific combinatorial games, like Wythoff's game and Chomp.

%%%%%%%%%%%%%%%%%%%%%%%%%%%%%%%%%%%%%%%%%%%%%%%%%%%%%%%%%%%%%%%%%%%%%%%%%
%%%%%%%%%%%%%%%%%%%%%%%%%%%%%%%%%%%%%%%%%%%%%%%%%%%%
%%%%%%%%%%%%%%%%%%%%%%%%%%%%%%%%%%%%%%%%%%%%%%%%%%%%%%%%%%%%%%%%%%%%%%%%%

%%%%%%%%%%%%%%%%%%%%%%%%%%%%%%%%%%%%%%%%%%%%%%%%%%%%%%%%%%%%%%%%%%%%%%%%%
\end{document}